\documentclass[trsc,nonblindrev]{template/informs3} 

\OneAndAHalfSpacedXI 
\usepackage{hyperref}
\usepackage{accents}

\usepackage{graphicx}
\usepackage{subcaption}

\usepackage{natbib}
 \bibpunct[, ]{(}{)}{,}{a}{}{,}%
 \def\bibfont{\small}%

 \usepackage{dsfont}
 \usepackage{diagbox}
 \usepackage{tikz}
 \usetikzlibrary{patterns}
 \usepackage{multicol}
 \usepackage{booktabs}
 \usepackage{xcolor,colortbl}
 \usepackage{rotating}

\usepackage{multirow}
\usepackage{pgfplots}
\pgfplotsset{compat=1.8}
\usepgfplotslibrary{statistics}

\definecolor{flowblue1}{HTML}{DEEBF7}
\definecolor{flowblue2}{HTML}{D6DCE5}
\definecolor{flowgray}{HTML}{D0CECE}
\usetikzlibrary{shapes.geometric}
\usetikzlibrary{matrix,calc,shapes}
\usetikzlibrary{intersections} 
\tikzset{
  treenode/.style = {shape=rectangle, line width = 0.25pt, 
                     draw = gray, anchor=center,
                     minimum width=10.25em, minimum height=3.5em, align=center,
                     inner sep=1.5ex},
  decision/.style = {treenode, diamond, inner sep=0pt},
  root/.style     = {treenode, font=\Large, bottom color=red!30},
  env/.style      = {treenode, font=\ttfamily\normalsize},
  finish/.style   = {root, bottom color=green!40},
  dummy/.style    = {circle,draw}
}

\TheoremsNumberedThrough      

\EquationsNumberedThrough    

\usepackage{wrapfig}
\usepackage{caption}
\usepackage{pdflscape}
\newcounter{alignnumber}
\renewcommand{\theequation}{\arabic{alignnumber}\alph{equation}}

\definecolor{Tblue}{HTML}{0065BD}
\definecolor{TUMgray}{HTML}{CCCCCC}
\definecolor{TUMdarkgray}{HTML}{808080}
\definecolor{TUMblue2}{HTML}{005293}
\definecolor{TUMblue3}{HTML}{003359}
\definecolor{TUMalblack}{HTML}{333333}
\definecolor{TUMorange}{HTML}{E37222}

\usepackage[ruled,vlined]{algorithm2e}
\SetAlgorithmName{Procedure}{Procedure}{List of Procedures}

\begin{document}


\RUNAUTHOR{Bauerhenne, Bard, and Kolisch}

\RUNTITLE{Robust Routing and Scheduling of Home Healthcare Workers}

\TITLE{Robust Routing and Scheduling of Home Healthcare Workers: A Nested Branch-and-Price Approach}

\ARTICLEAUTHORS{%
\AUTHOR{Carolin Bauerhenne$^\star$}
\AUTHOR{Jonathan Bard$^\dagger$}
\AUTHOR{Rainer Kolisch$^\diamond$}
}

\newpage
\TITLE{Robust Routing and Scheduling of Home Healthcare Workers: A Nested Branch-and-Price Approach}

\ABSTRACT{%
The global home healthcare market is growing rapidly due to aging populations, advancements in healthcare technology, and patient preference for home-based care. In this paper, we study the multi-day planning problem of simultaneously deciding patient acceptance, assignment, routing, and scheduling under uncertain travel and service times. Our approach ensures cardinality-constrained robustness with respect to timely patient care and the prevention of overtime.  
We take into account a wide range of criteria including patient time windows, caregiver availability and compatibility, a minimum time interval between two visits of a patient, the total number of required visits, continuity of care, and profit. We use a novel systematic modeling scheme that prioritizes health-related criteria as hard constraints and optimizes cost and preference-related criteria as part of the objective function. We present a mixed-integer linear program formulation, along with a nested branch-and-price technique. Results from a case study in Austin, Texas demonstrate that instances of realistic size can be solved to optimality within reasonable runtimes. The price of robustness primarily results from reduced patient load per caregiver. Interestingly, the criterion of geographical proximity appears to be of secondary priority when selecting new patients and assigning them to caregivers.  
}

\KEYWORDS{home healthcare; robust optimization; branch and price}

\maketitle
\section{Introduction} \label{intro} 
\label{introduction_HHC}

Home healthcare accounted for 3.5\% of national health expenditures in the United States in 2019 \citep{healthexpenditures}. 
The global home healthcare market, valued at USD 336.0 billion in 2021, is predicted to expand at a compound annual growth rate of 7.9\% from 2022 to 2030 \citep{marketgrowth}. This substantial growth is attributed to factors like aging populations, advancements in healthcare technology, and a rising preference for home-based care by patients and their families \citep{mckinsey}. Additionally, home healthcare is typically a more cost-effective option compared to hospital care, as it involves fewer resources \citep{singh2022comprehensive}. 

Routing and scheduling caregivers for home visits is a complex but essential task for operational smoothness in home healthcare. Firstly, the process involves carefully balancing the needs of patients, caregivers, and the agency, encompassing timely care delivery, favorable working conditions, and cost-efficiency. Secondly, the complexity arises from the high number of decisions, which make it difficult to find optimal solutions within reasonable computation times. Thirdly, the uncertainty in patient conditions and traffic further complicates the problem. 

In this paper, we study the problem of simultaneously deciding on patient selection, assignment, routing and scheduling over a multi-day planning horizon with uncertain travel and service times. We find robust schedules that ensure timely care and prevent overtime, even with a worst-case realization of service and travel times. The conservatism of the schedules can be controlled through the utilization of uncertainty budgets, as proposed by \cite{bertsimas2004price}. Our approach takes into account a wide range of optimization criteria including patient time windows, number of required visits, caregiver availability, compatibility, continuity of care, and profit.  
The main contributions of this study can be summarized as follows:

\paragraph{Robustness.} Building upon previous studies that primarily emphasized robustness concerning overtime, our research extends this focus to address the crucial element of ensuring timely patient care. The computational complexity increases as our approach requires meeting all patients' time windows in each worst-case scenario.

\paragraph{Health-outcome prioritizing modeling scheme.} We introduce a systematic modeling scheme for routing and scheduling in home healthcare. Previous studies have implicitly prioritized multiple criteria in various forms, as can be seen in the literature reviews by \cite{fikar2017home} and \cite{cisse2017or}. In contrast, we strictly differentiate between health-related criteria (modeled as hard constraints) and efficiency-related criteria (optimized in the objective function). This clear distinction ensures that all accepted patients receive the best possible treatment regarding factors influenced by operational planning.

\paragraph{Optimal solutions with branch and price.} 
We propose a nested Dantzig-Wolfe decomposition to tackle the complex structure of home healthcare routing and scheduling (HHCR\&S). Specifically, we decompose the problem into smaller subproblems for each caregiver, and subsequently into subproblems for each day. Based on this nested structure, we present the first branch-and-price method to solve a robust HHCR\&S problem. Notably, this method also represents one of the first branch-and-price approaches for multi-period HHCR\&S. 

\paragraph{Case study.} 
Using real data from a home healthcare agency in Dallas, Texas, our computational study demonstrates that the nested branch-and-price procedure yields optimal robust solutions for multi-day problem instances of realistic size within reasonable runtimes. While non-robust solutions are found within seconds for most instances, robustness increases runtimes to several minutes. Moreover, compared to parameter settings lacking uncertainty or time windows, introducing a medium level of robustness and one-hour time windows leads to an average profit loss of 14.20\%, primarily due to reduced patient loads per caregiver per day. Nonetheless, incorporating buffer time allows for low-cost robustification of some weekly schedules. Notably, our findings underscore that revenue takes precedence over geographic considerations in acceptance and assignment decisions.

The remainder of the paper is organized as follows. Section~\ref{literature_review} reviews related optimization problems in the HHCR\&S literature. Section~\ref{model} presents a robust mixed-integer linear model for the problem of weekly HHCR\&S. In Section~\ref{method}, we develop a branch-and-price algorithm to solve the problem efficiently also for larger problem instances. In Section~\ref{computational_study}, we conduct a computational study using real-world data from a home healthcare agency based in Dallas, Texas.

\section{Literature Review}\label{problemdescription} 
\label{literature_review}
For a comprehensive survey of the HHCR\rm{\&}S literature, we refer to \cite{cisse2017or, fikar2017home}, and \cite{di2021routing}. This review focuses on two streams that are most relevant to our research: In Section~\ref{Ereview}, we present an overview of exact solution methods for HHCR\rm{\&}S problems. In Section~\ref{Rreview}, we discuss robust optimization approaches in the presence of uncertain travel and/or service times. Our contribution to the literature lies in expanding the problem scope within each stream and innovatively combining both streams of research.

\subsection{Exact solution methods}
\label{Ereview}
To find optimal solutions for larger instances of HHCR\rm{\&}S problems, exact solution methods leverage the decomposable problem structure. Specifically, for a given assignment of patients to caregivers, the routing and scheduling subproblems for each caregiver can be solved more efficiently. Based on this observation, two well-known exact methods have been applied in the literature: Logic-based Benders decomposition (LBBD) and Branch and Price (B\rm{\&}P).

\paragraph{LBBD.}
HHCR\rm{\&}S can be elegantly solved using Logic-based Benders Decomposition. This approach involves dividing the problem into a master problem that assigns patients to caregivers, and subproblems that determine the routing and scheduling for each caregiver. 

The pioneering work by \cite{heching2016logic} introduces LBBD to HHCR\rm{\&}S. Specifically, they study a multi-period problem aiming to maximize the number of patients served, while considering time windows and continuity of care. They address the master problem through mixed-integer linear programming (MILP) and solve the subproblems by constraint programming. In subsequent work by \cite{heching2019logic}, they substantially enhance efficiency by introducing a time window subproblem relaxation, classical assignment relaxation, and multicommodity flow relaxation into the master problem; and introduce a variant of LBBD called branch and check.
\cite{grenouilleau2020new} propose a new algorithm to handle the subproblems and present a novel LBBD formulation featuring variables representing pre-generated visit patterns. While they do not aim for exact solutions, their matheuristic focuses on efficiency. In contrast, \cite{naderi2023novel} recently achieved an exact solution for a robust multi-period HHCR\rm{\&}S problem using LBBD. Their approach ensures robustness regarding overtime, but does not account for time windows. They minimize fixed caregiver costs and overtime costs. Notably, their LBBD approach differs from \cite{heching2019logic} due to the inclusion of a nurse eligibility feature. Moreover, their subproblems are optimization problems rather than feasibility problems. Consequently, they additionally developed custom Benders optimality cuts.

Notably, the inherent objective function of the LBBD master problem does not account for subproblem variables. For example, the travel cost is dependent on the routes between patients and cannot be solely determined by assigning patients to caregivers. Even though the master problem's objective function can be modified to incorporate objective functions based on subproblem variables, if the original problem's objective function aims to minimize travel or wage costs, the LBBD master problem may still not serve as the most effective guide to achieving an optimal overall solution. 

\paragraph{B \rm{\&}\it{P}.}
In the context of branch and price, the master problem selects either weekly or daily routes for each caregiver from a subset of all feasible solutions. The pricing problem, meanwhile, generates additional routes with better objective function values, if possible, and adds them to the pool of feasible solutions in the master problem. To the best of our knowledge, \cite{rasmussen2012home} are the first to employ a branch-and-price algorithm to a HHCR\rm{\&}S problem, maximizing the number of patients served for a single-period planning horizon. In contrast, \cite{trautsamwieser2014branch} address a multi-period context and minimize total working time. \cite{yuan2014home} and \cite{yuan2015branch} tackle the minimization of total cost for a single-period planning horizon, accounting for stochastic service times and late arrivals, using branch and price. \cite{liu2017mathematical} and \cite{liu2019branch} refine this approach to minimize travel costs and penalties for unvisited patient, with the former considering lunch breaks and the latter incorporating stochastic service and travel times. \cite{qiu2022exact} adds synchronization constraints, using a strengthened branch-and-price-and-cut algorithm. Similarly, \cite{tanoumand2021exact} concentrates on travel cost minimization with different vehicle types. \cite{alkaabneh2022multi} explore a multi-objective scenario, aiming to minimize cost and maximize patient-caregiver compatibility, combining branch and price with a two-stage metaheuristic. 

Recently, \cite{yin2023exact} tailored an efficient branch-and-price-and-cut algorithm for a single-period planning problem, considering dispatching and transport costs, as well as incompatibility costs, and skill and time window constraints. Their B\rm{\&}P algorithm involves a bounded bidirectional labelling algorithm and heuristics to solve the subproblems, as well as subset-row inequalities. These studies collectively highlight the efficiency of branch-and-price approaches in addressing the complexities of HHCR\rm{\&}S.

While some of these approaches incorporate stochastic travel and service times, none considers robustness with respect to unfavorable realizations. In related problems, such as the Robust Vehicle Routing Problem with Time Windows, prior research has demonstrated the effectiveness of branch and price \citep{pecin2017new, munari2019robust}.

Our contributions to exact solution methods for HHCR\rm{\&}S are twofold. First, we build upon a problem framework similar to \cite{heching2019logic}, extend it to robustness under uncertain service and travel times, and incorporate the objective of minimizing both travel and wage costs while simultaneously maximizing the number of accepted patients. As we consider minimizing travel cost a substantial objective, we employ B\rm{\&}P as our solution approach. Second, our contributions to the B\rm{\&}P literature encompass two key aspects: introducing robustness for travel and service times, and extending solutions to a multi-day planning horizon.

\subsection{Robust home healthcare schedules}
\label{Rreview} 
Disruptions such as traffic delays and changing patient needs pose recurring challenges to timely service delivery. While the prevailing literature in home healthcare focuses on deterministic routing and scheduling, recent years have witnessed a growing research interest in stochastic and robust planning. Given the direct impact of both timely care and overtime on patient health outcomes \citep{aiken2002hospital, bae2014assessing, hermida2016chronotherapy}, our study aims to ensure punctuality not only on average but also in worst-case scenarios. In the following, we present a comprehensive review of robust approaches for HHCR\rm{\&}S. Our examination aligns with the cardinality-constrained robustness framework (or $\Gamma$-robustness) introduced by \cite{bertsimas2004price}. For fuzzy approaches, we point to \cite{fathollahi2020bi} and \cite{shiri2023sustainable}. Regarding robust methods with discrete scenario sets, which are computationally easier to integrate, we refer to \cite{naji2017robust}, \cite{fathollahi2021multi}, \cite{fathollahi2022sustainable}, and \cite{shiri2023sustainable}.

Initial robust optimization studies addressed uncertain visit demand \citep{carello2014cardinality, carello2017handling, cappanera2017pattern}, which clearly is important for the assignment, but does not ensure timely care. Pioneering cardinality-constrained uncertainty in travel and service times, \cite{shi2019robust} tackle a single-period home healthcare scheduling problem, minimizing costs under time windows and skill requirements. They employed a recursive approach to linearize the adversarial problem. They find optimal solutions with a MILP solver and additional heuristic solutions with Simulated Annealing, Tabu Search and Variable Neighborhood Search. \cite{makboul2021robust} consider a similar single-period problem, but additionally take lunch breaks into account. They solve the problem using a MILP solver. \cite{shahnejat2021robust} extend the problem to multiple time windows, synchronization, and precedence constraints. They do not linearize the adversarial problem using the recursive approach, but enumerate all extreme points and solve the problem with metaheuristic algorithms. \cite{fathollahi2020robust} are the first to address a multi-period problem. They minimize total costs and carbon emissions, integrate robustness, and also solve the problem with metaheuristic algorithms. Finally, \cite{naderi2023novel} solve the first multi-period robust staffing, assignment, routing and scheduling problem exactly. Their goal is to minimize fixed costs and overtime costs. They employ dualization to linearize the adversarial problem of maximizing overtime, and do not consider time windows. They use pre-generated visit patterns and a novel logic-based Benders branching-decomposition algorithm for exact solutions.

This paper contributes to cardinality-constrained robust home healthcare routing and scheduling approaches as follows. Firstly, we ensure timely care within designated time windows on a multi-period planning horizon. Secondly, we conduct an analytical comparison of two linearization techniques for handling the robustness constraints. Thirdly, this study introduces a novel branch-and-price algorithm to find robust home healthcare routes and schedules, efficiently solving realistic problem instances to optimality.

\section{Robust Home Healthcare Routing and Scheduling}
\label{model}
\vspace{8pt}

\subsection{Problem statement}

Given a set of patients and a set of caregivers, we address the problem of simultaneously determining the acceptance of new patients, assigning patients to caregivers, and creating daily routes and schedules for a set of patients and a set of home health caregivers, over a multi-day planning horizon. We propose a new classification scheme for this multi-criteria problem whose objective is to optimize cost and preference considerations. All criteria related to health outcome, such as continuity of care, punctuality, and compatible assignments, are met as hard constraints to ensure the best possible health outcomes. We note that quality measures related to health outcomes are also the primary determinant for home healthcare ratings \citep{cms}. 
Specifically, our model takes into account the following criteria.

\subsubsection{Cost-related criteria.} 
\paragraph{Maximum profit.} To determine the total profit, we calculate the revenue generated from treating patients and subtract wage costs and travel costs.

\subsubsection{Health-outcome related criteria.} 
\paragraph{Patient time windows.} A majority of patients in home healthcare require assistance with activities of daily living, including bathing, dressing, toileting, transferring, eating, and administering medications \citep{caffrey2011home, osakwe2019activities}. Delays in these services can harm the patient's safety. Moreover, improper timing of certain medications, such as blood pressure medications or anti-diabetic medications, can result in serious health complications and hospitalizations \citep{hermida2016chronotherapy}. Therefore, we ensure that all patient time windows must be met even in a worst-case scenario of travel and service times.

\paragraph{Overtime.}  Overtime not only negatively impacts the physical and mental health of caregivers, but also the safety of patients and their health outcomes, see, e.g., \cite{aiken2002hospital, bae2014assessing}. Therefore, we ensure that caregivers have no overtime, even in a worst-case scenario of travel and service times.

\paragraph{Continuity of care.} 
Interpersonal continuity of care not only builds trust and increases patient satisfaction, it also contributes to improved health outcome and lower morbidity \citep{cabana2004does, chan2021effects}. This may be due, e.g.,  to a better understanding of the patient's medical history and treatment plan. We therefore enforce strict continuity of care, i.e., a patient is only accepted by the agency if he or she can be assigned to exactly one caregiver over the episode of care.

\paragraph{Compatible assignments.} 
Patient assignments must take into account the compatibility of the caregivers' skill levels, as well as other factors that may affect compatibility, such as gender and language.

\subsection{Model formulation}

\paragraph{Indices.} 
The planning horizon consists of multiple days, denoted by $ d \in \mathcal{D}$. The available caregivers are denoted by $k \in \mathcal{K}$. There are two types of patient requests: new patient requests for treatment than can be admitted or rejected, denoted by $i \in \mathcal{I}^{\mathcal{R}}$, and patients who must receive treatment, denoted by $i \in \mathcal{I}^{\mathcal{E}}$, including existing patients with ongoing treatments from the previous planning period. Patients in the set $\mathcal{I}^{\mathcal{E}}$ are pre-assigned to the same caregiver and days as in the previous planning period to ensure strict continuity of care. For each caregiver $k$ and day $d$, these patients compose the set $\mathcal{I}^{\mathcal{E}}_{k,d}$. All patients are grouped together in the set $\mathcal{I} := \mathcal{I}^{\mathcal{E}} \cup \mathcal{I}^{\mathcal{R}}$, and the total number of patients in $\mathcal{I}$ is $n$. We denote $N := \lbrace 1, \ldots, n \rbrace$ and $N_0 := \lbrace 0 \rbrace \cup N$. The depot, which is the starting and ending point for all caregivers, is indexed by $i \in \lbrace 0 , n+1 \rbrace$, respectively. 

\paragraph{Parameters.} Caregiver $k$ is available on day $d$ during a working time window of $\left[\underaccent{\bar}{W}_{ {  kd} }, \bar{W}_{ {  kd} }\right]$. Each patient $i$ has a time window $[\underaccent\bar{T}_{ {  id} },\bar{T}_{ {  id} }]$ on day~$d$. 
Moreover, patient $i$ requires exactly $v_i$ visits within the planning horizon, with at least $d_i$ days between any two visits. A compatibility matrix $C \in \lbrace0,1\rbrace^{n\times|\mathcal{K}|\times|\mathcal{D}|}$ indicates which patients are compatible with which caregivers and on which days. Specifically, caregiver $k$ can treat patient $i$ on day $d$ if and only if $C^i_{kd} = 1$. 
The revenue generated by caregiver $k$ when treating patient $i$ is denoted by $R_{ {  ki} }$. Furthermore, we take into account per time unit wage costs $C^W_{ {  k} }$ and travel costs $C^T_{ {  dij} }$ between two locations $i$ and $j$.

\paragraph{Uncertainty.}
We aim to find robust schedules in the presence of uncertain travel times $\tilde{t} \in \mathbb{R}_{+}^{(n+2) \times (n+2)}$ and uncertain service times $\tilde{p} \in \mathbb{R}_{+}^n$. 
For each patient $i \in \mathcal{I}$, we are given the expected service time $\bar{p}_{ {  i} }$ and the maximum positive deviation $\hat{p}_{ {  i} }$. Similarly, for each transfer from location $i$ to location $j$, we are given the expected travel time $\bar{t}_{ {  ij} }$ and the maximum positive deviation $\hat{t}_{ {  ij} }$. The uncertainty budget for travel times is denoted by $\Gamma^t \in \mathbb{N}_{0}$, and the uncertainty budget for service times is denoted by $\Gamma^p \in \mathbb{N}_{0}$. Following the polyhedral uncertainty sets proposed in \cite{bertsimas2004price}, with the worst-case being achieved at the right-hand side of the uncertainty range, the uncertainty sets for travel and service times are defined as 
$\mathcal{U}^p := \left\lbrace p \in \mathbb{R}_{+}^n : p_{ {  i} } = \bar{p}_{ {  i} } + \gamma_{ {  i} }^p \hat{p}_{ {  i} }, \sum_{i=1}^n \gamma_{ {  i} }^p \leq \Gamma^p, 0 \leq \gamma_{ {  i} }^p \leq 1 \,\forall i \in \lbrace 1, \ldots, n \rbrace \right\rbrace$
and $\mathcal{U}^t := \left\lbrace t \in \mathbb{R}_{+}^{(n+2) \times (n+2)} : t_{ {  ij} } = \bar{t}_{ {  ij} } + \gamma_{ {  ij} }^t \hat{t}_{ {  ij} }, \sum_{i=0}^{n+1} \sum_{i=0}^{n+1}  \gamma_{ {  ij} }^t \leq \Gamma^t, 0 \leq \gamma_{ {  ij} }^t \leq 1 \,\forall i,j \in \lbrace 0, \ldots, n+1 \rbrace \right\rbrace$
respectively. We refer to each element in the uncertainty sets of travel and service times as a travel or service time scenario, respectively. In particular, all uncertainty is assumed to be uncorrelated.
The uncertainty budgets control the level of conservatism, and will be shown to correspond to the daily number of patients (transfers) per caregiver that require the maximum absolute deviation from their expected service (travel) time.

\paragraph{Variables.} 
We use binary variables $x_{ {  ij} }^{ {  kd} }$ to determine whether caregiver $k \in \mathcal{K}$ on day $d \in \mathcal{D}$ treats patient $j \in \mathcal{I} \cup \lbrace 0, n + 1 \rbrace$ immediately after patient $i\in \mathcal{I} \cup \lbrace 0, n + 1 \rbrace$. We use binary variables $\delta_{ki} \in \lbrace 0, 1 \rbrace$ to determine whether patient $i \in \mathcal{I}$ is assigned to caregiver $k \in \mathcal{K}$.
Moreover, we find the scenario-dependent start time of patient $i \in \mathcal{I}$ on day $d$ as:
\[
s_{i,d}(x,p,t) :=   \min \left\{
\begin{array}{l}
    \underaccent\bar{T}_{ {  id} },  \\
    \max \left\lbrace \begin{array}{l}
    \sum\limits_{k \in \mathcal{K}} \left( \underaccent{\bar}{W}_{ {  kd} } + t_{ {  0i} }\right) x_{ {  0i} }^{ {  kd} },  \\
     \max\limits_{j \in \mathcal{I}} \left\lbrace s_{ {  j} } + (p_{ {  j} } + t_{ {  ji} })\sum\limits_{k \in \mathcal{K}}x_{ {  ji} }^{ {  kd} } + (\underaccent\bar{T}_{ {  id} }  - \bar{T}_{ {  jd} })(1 - \sum\limits_{k \in \mathcal{K}} x_{ {  ji} }^{ {  kd} }) \right\rbrace \\
\end{array}\right\}
\end{array}\right\}. 
\]

The first term of the minimization function indicates that patient $i$ is available for treatment starting with the lower bound $\underaccent\bar{T}_{ {  id} }$ of his or her time window. The second term in the minimization function states that caregiver $k$ starts working at time $ \underaccent{\bar}{W}_{ {  kd} }$, and is available to treat patient $i$ only after completing treatment of the previous patient and traveling to patient $i$. 

Similarly, we define the scenario-dependent return time of caregiver $k$ to the depot on day $d$ as follows:
\[
s^{n+1}_{kd}(x,p,t) := \min \left\{
\begin{array}{l}
    \underaccent{\bar}{W}_{ {  kd} },  \\
    \max\limits_{i \in \mathcal{I}} \left\lbrace s_{ {  i, d} }(x,p,t) + (p_{ {  i} } + t_{i,n+1})x_{i,n+1}^{kd} + (\underaccent{\bar}{W}_{ {  kd} }  - \bar{T}_{ {  id} })(1 -  x_{i,n+1}^{kd}) \right\rbrace \\
\end{array}\right\}.
\]

The first term of the minimization function ensures that  caregiver $k$ cannot return to the depot before starting his or her shift. The second term in the minimization function states that the caregiver returns to the depot only after completing the treatment of the last patient and traveling to the depot. 
 
We now consider the following Robust Home Healthcare Routing and Scheduling Problem (R-HHCR\rm{\&}S). \label{R-HHCRS}

\begin{align}     \stepcounter{alignnumber}     \setcounter{equation}{0}  
    \max\, & \multicolumn{2}{l}{$\sum\limits_{k\in\mathcal{K}} \left( \sum\limits_{\substack{i \in \mathcal{I}}} R_{ki}v_i\delta_{ki} -   \sum\limits_{d\in\mathcal{D}} \sum\limits_{\substack{i, j \in \mathcal{I}}}C^T_{dij}  x_{ {  ij} }^{ {  kd} } -   C^W_k  \sum\limits_{d\in\mathcal{D}} \left({s}^{n+1}_{kd}-\underaccent{\bar}{W}_{ {  kd} }\right)
        \right) $} \label{model1obj} \\[3pt]
     s.t. \, 
    & 
    \max\limits_{\scriptsize \substack{p \in \mathcal{U}^p \\  t \in \mathcal{U}^t}} s_{id}(x,p,t) \leq \bar{T}_{ {  id} } & \forall d \in \mathcal{D},i \in \mathcal{I} \label{model1con1} \\[3pt]
    & 
     \max\limits_{\scriptsize \substack{p \in \mathcal{U}^p \\  t \in \mathcal{U}^t}}  {s}^{n+1}_{kd}(x,p,t) =: {s}^{n+1}_{kd} \leq \bar{W}_{ {  kd} } & \forall k \in \mathcal{K} , d \in \mathcal{D} \label{model1con2} \\[3pt]
    &  \sum_{d\in\mathcal{D}} \sum_{j\in\mathcal{I}} x_{ {  ij} }^{ {  kd} } = v_i \delta_{ki} & \forall k \in \mathcal{K},i\in\mathcal{I} \label{model1con3} \\
        &      \sum_{k \in \mathcal{K}} \delta_{ki} \leq 1   & \forall i\in\mathcal{I} \label{model1con4} \\
                &      \sum_{i \in N_0} x_{ji}^{kd} = 1   & \forall k\in\mathcal{K} ,i\in\mathcal{I}^{E|_{k,d}} \label{model1con5} \\[3pt]
                                  &  \sum_{j\in\mathcal{I}} x_{ {  ij} }^{ {  kd} } \leq C_i^{kd}  & \forall k\in\mathcal{K},d\in\mathcal{D},i\in\mathcal{I}^{\mathcal{R}} \label{model1con6} \\
    & \sum_{k\in\mathcal{K}} \sum_{j\in\mathcal{I}} x_{ {  ij} }^{ {  kd} } + \sum_{k\in\mathcal{K}} \sum_{\mu \in [ d_i ]_0}  \sum_{j\in\mathcal{I}} x_{ij}^{kd+\mu} \leq 1 & \forall d\in\mathcal{D},i\in\mathcal{I}^{\mathcal{R}} \label{model1con7}  \\[3pt]
    & \sum_{j\in\mathcal{J}} x_{ {  ji} }^{ {  kd} }  -\sum_{j\in\mathcal{J}} x_{ {  ij} }^{ {  kd} } = 0, \quad \sum_{i\in\mathcal{J}} x_{ii}^{kd} = 0 & \forall k\in\mathcal{K},d\in\mathcal{D},i\in\mathcal{I} \label{model1con8} \\
    & \sum_{i \in [n+1]} x_{0,i}^{kd} = 1, \quad \sum_{i\in N_0} x_{i,n+1}^{kd} = 1 & \forall k\in\mathcal{K},d\in\mathcal{D} \label{model1con9} \\[3pt]
   & x_{ {  ij} }^{ {  kd} }, \delta_{ki} \in \lbrace 0,1 \rbrace & \forall k\in\mathcal{K},d\in\mathcal{D},i, j\in [n+1]_0 \label{model1con10}
\end{align}

The objective function~\hyperref[R-HHCRSP]{(\ref{model1obj})} maximizes the total profit, calculated by subtracting travel and wage costs from revenue. Constraints~\hyperref[R-HHCRSP]{(\ref{model1con1})} ensure robust schedules by guaranteeing that patient time windows are met even in a worst-case scenario of travel and service times. Similarly, constraints~\hyperref[R-HHCRSP]{(\ref{model1con2})} create robust schedules by ensuring that the caregivers' work hours are not exceeded even in a worst-case scenario of travel and service times. 
We call the maximization term within constraints~\hyperref[R-HHCRSP]{(\ref{model1con1})} and~\hyperref[R-HHCRSP]{(\ref{model1con2})} the {\it adversarial problems}. Due to their nonlinear nature, (\ref{model1obj}) -- (\ref{model1con10}) is not amenable to MILP solvers as it is currently formulated.

Constraints~\hyperref[R-HHCRSP]{(\ref{model1con3})} ensure that each accepted patient is visited the exact number of required times. By assigning each patient to at most one caregiver,  constraints~\hyperref[R-HHCRSP]{(\ref{model1con4})} enforce strict continuity of care. Constraints~\hyperref[R-HHCRSP]{(\ref{model1con5})} ensure that all patients in the set $\mathcal{I}^{\mathcal{E}}$ are visited by their pre-assigned caregiver on the designated day. Notably, adjusting the time windows allows for the straightforward fixing of specific visits to particular times. Constraints~\hyperref[R-HHCRSP]{(\ref{model1con6})} specify individual patient-caregiver compatibility.
Constraints~\hyperref[R-HHCRSP]{(\ref{model1con7})} enforce a minimum number of days between any two treatments for a patient.

Constraints~\hyperref[R-HHCRSP]{(\ref{model1con8})} are the well-known flow-conservation constraints. Constraints~\hyperref[R-HHCRSP]{(\ref{model1con9})} ensure that each caregiver starts and ends at a central depot. Together, the last two constraints ensure feasible routes, with constraints~\hyperref[R-HHCRSP]{(\ref{model1con6})} eliminating subtours. 

The model contains $\mathcal{O}\left(\mathcal{K}|\cdot |\mathcal{D}| \cdot n^2 \right) $ binary variables, $\mathcal{O} \left( |\mathcal{K}|\cdot |\mathcal{D}| \cdot n \right) $ continuous variables and $\mathcal{O}(n\cdot|\mathcal{K}|\cdot |\mathcal{D}|)$ constraints. 
After linearizing constraints~\hyperref[R-HHCRSP]{(\ref{model1con1})} and ~\hyperref[R-HHCRSP]{(\ref{model1con2})} it becomes a MILP.

\subsection{Adversarial problems}\label{adv} 

The \emph{dualization scheme} proposed by \cite{bertsimas2004price} has classically been used for linearizing adversarial problems in cardinality-constrained robust optimization. However, due to the nested definition of service start times, its application to our adversarial problems is not straightforward.
To address this issue, various alternative approaches have been proposed in the literature, such as by \cite{agra2013robust} and by \cite{lee2012robust}. 
One of the most efficient approaches currently available is the \emph{recursive scheme} introduced by \cite{munari2019robust} for a robust vehicle routing problem with time windows.

In the following, we compare the effectiveness of the classical dualization scheme, based on a novel reformulation of service start times, with the effectiveness of the recursive scheme for linearizing constraints~\hyperref[R-HHCRSP]{(\ref{model1con1})} and \hyperref[R-HHCRSP]{(\ref{model1con2})}. While the novel reformulation provides a deeper understanding of robust start times and how to choose the uncertainty budgets, our theoretical analysis demonstrates that the recursive scheme remains the more efficient method for linearizing the adversarial problems in \hyperref[R-HHCRS]{R-HHCR\rm{\&}S}.

\subsubsection{Dualization scheme.} 

To apply the classical dualization scheme as proposed by \cite{bertsimas2004price}, we reformulate the service start times. In particular, we replace the binary variable $x$ with a new binary variable $z_{j'j}^{k, d, i} \in \lbrace0,1\rbrace$, which additionally indicates the sequence of patients treated by caregiver $k$ on day $d$. Specifically, $z_{j'j}^{kdi}$ takes a value of 1 if and only if the $i$th patient treated by caregiver $k$ on day $d$ is patient $j$, and the patient treated before $j$ was patient $j'$; 0 otherwise. 
We now shift the focus from the start time $s_{id}$ of patient $i$ to the start time $s_{kd}^i$ of the $i$th patient treated by caregiver $k$ on day $d$. This start time can be calculated as follows (see Appendix~\ref{dual_scheme}): 

\begin{align}     \stepcounter{alignnumber}     \setcounter{equation}{0} 
\min \, & s_{kd}^{i} \label{reform} \\
   \text{s.t.} \,& s_{kd}^{i} \geq \underaccent{\bar}{W}_{ {  kd} }  + \max_{\substack{p \in \mathcal{U}^p \\t \in \mathcal{U}^t}} \left( \sum_{i' = 1}^{i - 1} \sum_{j'=0}^n\sum_{j=1}^n ( t_{j'j} + p_{ {  j} })z_{j'j}^{k,d,i'} + \sum_{j'=0}^n \sum_{j=1}^n t_{j'j}z_{j'j}^{k,d,i} \right) &  \label{max1} \\
   &  s_{kd}^{i} \geq  \sum_{j'=0}^n \sum_{j=1}^n \underaccent{\bar}{T}_{ {  jd} } z_{j'j}^{k,d,i^*}   + \max_{\substack{p \in \mathcal{U}^p \\t \in \mathcal{U}^t}} \sum_{i' = i^*+1}^{i} \sum_{j'=0}^n\sum_{j=1}^n (p_{ {  j'} } + t_{j'j})z_{j'j}^{k,d,i'} &  \forall i^* \in [i]  \label{max3}
\end{align}

Model (\ref{reform}) -- (\ref{max3}) ensures that the start time of each patient equals the maximum of the earliest availability of all previous visits (including the depot) plus a worst-case scenario of service and travel times to the current patient. This formulations allows us to now apply the dualization scheme (see Appendix~\ref{dual_scheme} for the derivation). 

Let $n_{max}$ denote the maximum number of patients treated by any caregiver on any day. The linearization of the robustness constraints~\hyperref[R-HHCRSP]{(\ref{model1con1})} and~\hyperref[R-HHCRSP]{(\ref{model1con2})} using the proposed dualization approach requires in total $n^2\cdot n_{max}\cdot|\mathcal{K}|\cdot|\mathcal{D}|$ binary variables, $\mathcal{O}\left(n^2\cdot \left(n_{max}\right)^2\cdot|\mathcal{K}|\cdot|\mathcal{D}| \right)$ continuous variables and $\mathcal{O}\left(n^2\cdot|\mathcal{K}|\cdot|\mathcal{D}|\cdot \left(n_{max}\right)^2 \right)$ new constraints. 
We note that if the uncertainty budget for travel time $\Gamma^t \geq 3$, then the proposed dualization scheme is more efficient compared to the linearization approach given by \cite{agra2013robust}.

\subsubsection{Recursive scheme.} \label{recurse}
Next, we apply the recursive scheme proposed by \cite{munari2019robust} to linearize the problem. Let $s_{(i, \gamma^p, \gamma^t)}^{d}$ denote the service start time of patient $i$ on day $d$, where $\gamma^p \in \lbrace 0, \ldots, \Gamma^p \rbrace$ and $\gamma^t \in \lbrace 0, \ldots, \Gamma^t \rbrace$ denote the units of maximum positive deviation that have already been consumed from the uncertainty budgets. Then we can replace the robustness constraints~\hyperref[R-HHCRSP]{(\ref{model1con1})} and~\hyperref[R-HHCRSP]{(\ref{model1con2})} by 
\begin{align}     \stepcounter{alignnumber}     \setcounter{equation}{0}  \label{rec1}
    s_{i, \gamma^p, \gamma^t}^{d} & \in [ \underaccent\bar{T}_{ {  id} }, \bar{T}_{ {  id} }] \\ & \qquad \forall d\in\mathcal{D},i \in \mathcal{I} ,\gamma^p \in [\Gamma^p]_0 ,\gamma^t \in [\Gamma^t]_0 \\
    s_{i, 0, \zeta^t}^{d} &  \geq  \sum_{k \in \mathcal{K}} \left(\underaccent{\bar}{W}_{ {  kd} }  + \bar{t}_{ {  0i} } + \zeta^t\hat{t}_{ {  0i} }\right)  x_{ {  0i} }^{ {  kd} }  \\ &  \qquad\forall d\in\mathcal{D},i\in\mathcal{I}\,  \forall  \zeta^t \in \lbrace 0, 1 \rbrace \\
    s_{i, \gamma^p + \zeta^p, \gamma^t + \zeta^t}^{d} &  \geq  s_{j, \gamma^p , \gamma^t}^{d} + (\bar{p}_j + \zeta^p\hat{p}_j + \bar{t}_{ {  ji} } + \zeta^t\hat{t}_{ {  ji} }) \sum_{k \in \mathcal{K}} x_{ {  ji} }^{ {  kd} } + ( \underaccent\bar{T}_{ {  id} } - \bar{T}_{ {  jd} })(1- \sum_{k \in \mathcal{K}} x_{ {  ji} }^{ {  kd} })  \\ &   \qquad \forall i,j\in\mathcal{I}  ,\gamma^p \in [\Gamma^p]_0 ,d\in\mathcal{D},\gamma^t \in [\Gamma^t]_0,\zeta^p, \zeta^t \in \lbrace 0, 1 \rbrace: \\ &  \qquad  \gamma^p + \zeta^p \leq \Gamma^p, \gamma^t + \zeta^t \leq \Gamma^t
\end{align}
and similarly, for avoiding overtime, 
\begin{align}     \stepcounter{alignnumber}     \setcounter{equation}{0} 
     {s}^{n+1}_{kd} & \in [\underaccent{\bar}{W}_{ {  kd} } , \bar{W}_{ {  kd} }]\\
     &  \qquad  \forall k\in\mathcal{K},d\in\mathcal{D} \\
    {s}^{n+1}_{kd} &  \geq  s_{i, \Gamma^p -\zeta^p, \Gamma^t -\zeta^t}^{d} + (\bar{p}_{ {  i} } + \zeta^p\hat{p}_{ {  i} } + \bar{t}_{i,n+1} + \zeta^t\hat{t}_{i,n+1})  x_{i,n+1}^{kd} + ( \underaccent{\bar}{W}_{ {  kd} } - \bar{T}_{ {  id} })(1- x_{i,n+1}^{kd}) \\ &  \qquad  \forall k\in\mathcal{K},d\in\mathcal{D},i\in\mathcal{I},\zeta^p, \zeta^t \in \lbrace 0, 1 \rbrace \label{rec2}
\end{align}
This leads to $\mathcal{O}\left( n\cdot|\mathcal{D}|\cdot|\Gamma^p|\cdot|\Gamma^t| + |\mathcal{K}|\cdot|\mathcal{D}| \right)$ new continuous variables, $\mathcal{O}\left(n^2\cdot|\mathcal{D}|\cdot|\Gamma^p|\cdot|\Gamma^t| \right)$ new constraints, but no new binary variables.
In contrast to the dualization scheme, the number of variables in the recursive scheme depends on the uncertainty budget. 

\subsubsection{Comparison.} \label{maxbudget} Based on the number of variables and constraints used in each approach, we finally compare the effectiveness of the dualization scheme and the recursive scheme to linearize our robustness constraints. 

Firstly, the recursive scheme stands out as it requires significantly fewer binary variables, an order of $n_{max}$ times less than the dualization scheme. Secondly, if $|\mathcal{K}|\cdot\left(n_{max}\right)^2 > |\Gamma^p|\cdot|\Gamma^t|$, then the recursive scheme also requires fewer continuous variables and constraints.

Despite this, examining the reformulation of start times in (\ref{reform}) -- (\ref{max3}) reveals that the uncertainty budgets limit the amount of uncertainty faced by each caregiver per day, and should be set accordingly. Therefore, w.l.o.g., we set the uncertainty budgets such that $ \Gamma^p\leq n_{max}$ and $\Gamma^t \leq n_{max} + 2$, a choice that further enhances the efficiency of the recursive scheme. It is worth noting that if the budgets exceed the maximum number of patients visited and paths traveled by a caregiver on a given day, the problem simplifies to a deterministic problem in which it is assumed that all patient visits and traveled paths require the maximum positive deviation in addition to the expected time. 
 
In conclusion, the recursive scheme proves more efficient in linearizing the robustness constraints in \hyperref[R-HHCRS]{R-HHCR\rm{\&}S}. As a result, we replace \hyperref[R-HHCRSP]{(\ref{model1con1})} -- \hyperref[R-HHCRSP]{(\ref{model1con2})} with the linear inequalities (\ref{rec1}) -- (\ref{rec2}).

\section{Solution Method}
\label{method}

To solve \hyperref[R-HHCRS]{R-HHCR\rm{\&}S} for larger instances, we exploit its decomposable problem structure. Specifically, 
\begin{figure}[b] \label{nest}
    \begin{tikzpicture}
    \node[rectangle, draw, minimum width = 0.95cm, minimum height = 0.35cm] (r) at (0.6,0.1) {};
	\node[rectangle, draw] (r) at (0.3,0.1) {};
	\node[rectangle, draw] (r) at (0.6,0.1) {};
	\node[rectangle, draw] (r) at (0.9,0.1) {};
	\node[rectangle, draw, minimum width = 0.95cm, minimum height = 0.35cm] (r) at (1.6,0.1) {};
	\node[rectangle, draw] (r) at (1.3,0.1) {};
	\node[rectangle, draw] (r) at (1.6,0.1) {};
	\node[rectangle, draw] (r) at (1.9,0.1) {};
	\node[rectangle, draw, minimum width = 0.95cm, minimum height = 0.35cm] (r) at (2.6,0.1) {};
	\node[rectangle, draw] (r) at (2.3,0.1) {};
	\node[rectangle, draw] (r) at (2.6,0.1) {};
	\node[rectangle, draw] (r) at (2.9,0.1) {};
    \node[rectangle, draw, minimum width = 0.95cm, minimum height = 1.25cm] (r) at (0.6,-0.75) {};
	\node[rectangle, draw] (r) at (0.3,-0.3) {};
	\node[rectangle, draw] (r) at (0.6,-0.3) {};
	\node[rectangle, draw] (r) at (0.9,-0.3) {};
	\node[rectangle, draw] (r) at (0.3,-0.6) {};
	\node[rectangle, draw] (r) at (0.6,-0.9) {};
	\node[rectangle, draw] (r) at (0.9,-1.2) {};
    \node[rectangle, draw, minimum width = 0.95cm, minimum height = 1.25cm] (r) at (1.6,-2) {};
	\node[rectangle, draw] (r) at (1.3,-1.55) {};
	\node[rectangle, draw] (r) at (1.6,-1.55) {};
	\node[rectangle, draw] (r) at (1.9,-1.55) {};
	\node[rectangle, draw] (r) at (1.3,-1.85) {};
	\node[rectangle, draw] (r) at (1.6,-2.15) {};
	\node[rectangle, draw] (r) at (1.9,-2.45) {};
	\node[rectangle, draw, minimum width = 0.95cm, minimum height = 1.25cm] (r) at (2.6,-3.25) {};
	\node[rectangle, draw] (r) at (2.3,-2.8) {};
	\node[rectangle, draw] (r) at (2.6,-2.8) {};
	\node[rectangle, draw] (r) at (2.9,-2.8) {};
	\node[rectangle, draw] (r) at (2.3,-3.1) {};
	\node[rectangle, draw] (r) at (2.6,-3.4) {};
	\node[rectangle, draw] (r) at (2.9,-3.7) {};
    \end{tikzpicture}
    \caption{Two-fold decomposition of the problem structure. \textmd{The small boxes represent daily routing and scheduling problems for each caregiver. The daily solutions for each caregiver are combined to create multi-day solutions, represented by the large boxes. The corresponding combination process ensures that the required number of treatments is met, and the minimum number of days between two visits is satisfied. The multi-day plans for each caregiver are further combined to ensure strict continuity of care and minimize overall cost. }
}
    \label{fig:problem_structure}
\end{figure}
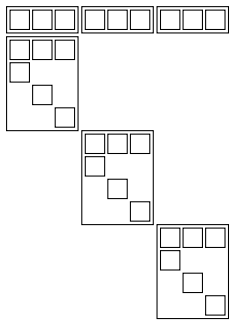
\begin{figure}[t] 
    \centering 
    \begin{tikzpicture}[remember picture, every node/.style={font=\footnotesize}]
    \tikzset{
        line/.style = {draw},
        comment/.style = {rectangle, draw, text centered, 
                rounded corners, minimum height=2em,fill=white},
        comment border/.style = {comment, rounded corners=0pt, 
                fill=gray!40, inner sep=4mm, minimum height=2em+4mm},
        terminator/.style = {shape=rounded rectangle, draw, inner sep=2mm, minimum width=12em, minimum height=3.25em, line width = 0.25pt},
}
    \filldraw[fill = gray!5!white, draw = gray!15!white] (-1.25, 0.6125) -- (15,0.6125) -- (15, 4.5) -- (-1.25, 4.5) -- cycle;
    \filldraw[fill = gray!5!white, draw = gray!15!white] (-1.25, .375) -- (15,.375) -- (15, -3.375) -- (-1.25, -3.375) -- cycle;
    \filldraw[fill = gray!5!white, draw = gray!15!white] (-1.25, -7.625) -- (15,-7.625) -- (15, -3.625) -- (-1.25, -3.625) -- cycle;
    \node at (-.75, 2.5) [rotate=90, font=\normalfont] {\color{gray} \bf Level 1};
    \node at (-.75, -1.5) [rotate=90, font=\normalfont] {\color{gray} \bf Level 2};
     \node at (-.75, -5.5) [rotate=90, font=\normalfont] {\color{gray} \bf Level 3};
     \node [align=center, terminator, draw = gray, fill = white] at (2.125,5.5) {Initialization};
     \draw[->, line width = 0.25pt, gray] (2.125,4.875) -- (2.125,3.22);
    \node [align=center, treenode, draw = gray, fill = white] at (2.125,2.5) {Relaxed master problem \\ \hyperref[MP-1]{(5)}};
    \draw[->, line width = 0.25pt, gray] (2.125,1.825) -- (2.125,-0.78);
    \node at (2, .5) [fill=white, inner sep = 1.5pt]{\color{gray} for each caregiver};
     \node [align=center, treenode, draw = gray, fill = white] at (2.125,-1.5) {Relaxed master problem \\ \hyperref[MP2-first]{(7)}};
     \draw[->, line width = 0.25pt, gray] (2.125,-2.195) -- (2.125,-3.975);
      \node at (2.125, -3.5) [fill=white, inner sep = 1.5pt]{\color{gray} for each day};
     \node [align=center, treenode, draw = gray, fill = white] at (2.125,-5.5) {Branch and bound \\ acceptance decision by \\ \hyperref[knapsack]{multi-Knapsack} problem \\ with \hyperref[TSP]{R-TSPTW} \\ subproblems};
     \draw[->, line width = 0.25pt, gray] (4.125,-5.5) -- (4.84,-5.5);
     \node [align=center, terminator, draw = gray!20!white, fill = white, rounded corners] at (7,-5.5) {Daily caregiver's  plan};
     \draw[->, line width = 0.25pt, gray] (7,-4.875) -- (7,-4.425);
     \node [align=center, diamond, draw = gray, fill = white, aspect=3] at (7,-3.5) {Cost-reducing \\ daily plan?};
     \draw[line width = 0.25pt, gray] (3.375,-3.5) -- (4.35,-3.5);
     \draw[->, line width = 0.25pt, gray] (3.375,-3.5) -- (3.375,-2.25);
     \draw[->, line width = 0.25pt, gray] (7,-2.65) -- (7,-2.25);
      \node at (4.35, -3.5)[fill=white, inner sep = 0.5pt]{\color{gray} Y};
      \node at (4.1, -2.65) [fill=gray!5!white, inner sep = 1.5pt]{\color{gray} Add daily plan};
      \node at (7, -2.65) [fill=gray!5!white, inner sep = 1.5pt]{\color{gray} N};
     \node [align=center, treenode, draw = gray, fill = white] at (7,-1.5) {Find integral solution via \\ branch and bound};
     \draw[->, line width = 0.25pt, gray] (4.95,-1.5) -- (4.19,-1.5);
     \draw[->, line width = 0.25pt, gray] (9.05,-1.5) -- (9.85,-1.5);
      \node [align=center, terminator, draw = gray!20!white, fill = white, rounded corners] at (12,-1.5) {Weekly caregiver's  plan};
      \draw[->, line width = 0.25pt, gray] (12,-.85) -- (12,-.5);
      \node [align=center, diamond, draw = gray, fill = white, aspect=3] at (12,0.5) {Cost-reducing \\ caregiver's plan?};
       \draw[line width = 0.25pt, gray] (3.375,0.5) -- (9.1,0.5);
     \draw[->, line width = 0.25pt, gray] (3.375,0.5) -- (3.375,1.79);
      \draw[->, line width = 0.25pt, gray] (12,1.48) -- (12,1.79);
      \node at (9.1, 0.5) [fill=white, inner sep = 0.5pt] {\color{gray} Y};
      \node at (4.5, 1.48) [fill=gray!5!white, inner sep = 1.5pt]{\color{gray} Add caregiver's plan};
      \node at (12, 1.48) [fill=gray!5!white, inner sep = 1.5pt]{\color{gray} N};
      \node [align=center, treenode, draw = gray, fill = white] at (12,2.5) {Find integral solution via \\ branch and bound};
       \draw[->, line width = 0.25pt, gray] (9.95,2.5) -- (4.19,2.5);
    \draw[->, line width = 0.25pt, gray] (12,3.17) -- (12,4.825);
     \node [align=center, terminator, draw = gray, fill = white,] at (12,5.5) {Optimal \\ \hyperref[R-HHCRS]{R-HHCR\rm{\&}S} plan};
\end{tikzpicture}
\centering
    \caption{Flowchart of the nested branch-and-price procedure}
    \label{flowchart}
\end{figure}
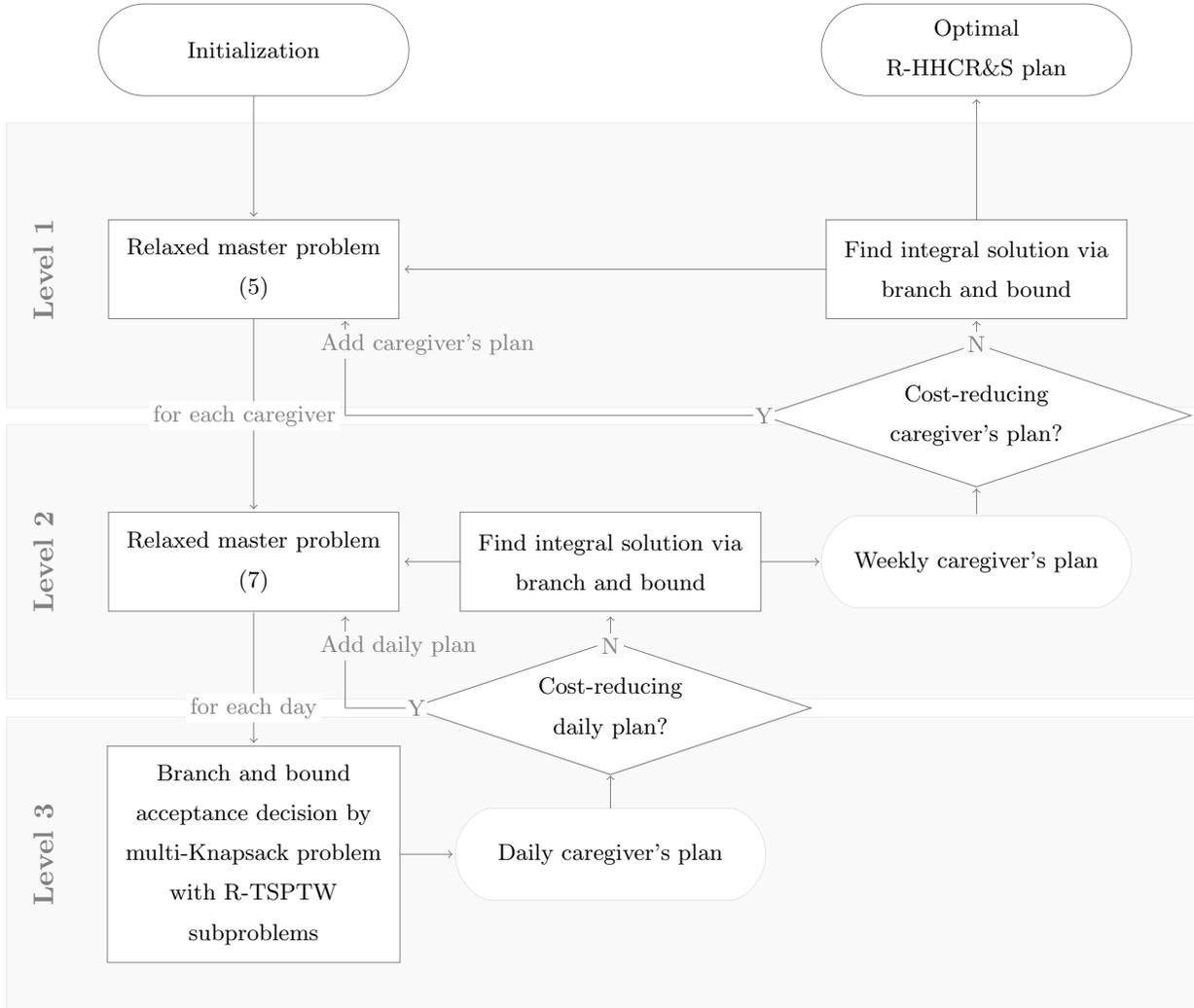
the problem can be decomposed into subproblems for each caregiver, and each caregiver's subproblem can be further decomposed into daily routing and scheduling subproblems, as depicted in Figure~\ref{fig:problem_structure}. Our nested branch-and-price procedure, outlined in Figure~\ref{flowchart}, leverages this decomposition to recursively compute solutions across three levels.

At the first level, the restricted master problem combines given plans for all caregivers into a comprehensive weekly plan for the agency. Integral solutions for the selection of caregiver plans are ensured via a branch-and-bound procedure. In the second level, new weekly plans for each caregiver are generated, terminating either with a caregiver's weekly plan that potentially improves the agency's overall plan or with the proof that none exists. Weekly caregiver schedules are generated using the second-level relaxed master problem, which combines daily plans into a weekly schedule. Integral solutions for daily route selection are again ensured through a branch-and-bound procedure. This is accomplished using a branch-and-bound procedure with robust traveling salesman subproblems incorporating time windows, the acceptance decision being guided by a multi-Knapsack problem.

Overall, our nested branch-and-price procedure efficiently explores promising daily routes and weekly plans, combining them into a comprehensive agency plan. Notably, while nested decomposition approaches have been employed in related problems \citep{vanderbeck2001nested, hennig2012nested, dohn2013branch, tilk2019nested, zheng2021stochastic}, to the best of our knowledge we are the first to apply this method to an HHCR\rm{\&}S problem. 

\subsection{Nested Dantzig-Wolfe decomposition}

\subsubsection{First-level decomposition.} 
We first decompose \hyperref[R-HHCRS]{R-HHCR\rm{\&}S} into a master problem and $|\mathcal{K}|$ pricing problems. Each pricing problem, associated with a fixed caregiver $k \in \mathcal{K}$, determines a weekly plan $P$ that encodes a schedule and route for each day. We note that caregivers are heterogeneous in compatibility to patients, work times, wage cost, and revenue. The master problem then combines the weekly plans generated by the pricing problems in order to maximize profit, while ensuring that care is provided to all accepted patients.

\paragraph{Notation.} We use a binary variable $\lambda_{k,P}$ to decide whether caregiver $k$ is assigned to the weekly plan~$P$. The set of all feasible weekly plans for caregiver $k$ is denoted as $\mathcal{P}_k$. We use the binary parameter $a^i_{k,P}$ to indicate whether patient $i$ is being treated by caregiver $k$ in plan $P$. The revenue generated by a plan $P$ executed by caregiver $k$ is denoted by $Rev_{k,P}$, and similarly the travel cost $C^T_{k,P}$ and wage cost $C^W_{k,P}$.

\paragraph{Master problem.}
The master problem is a variant of the set partitioning problem, in which elements are either accepted or rejected, and the accepted elements are partitioned into $|\mathcal{K}|$ subsets. Similar to other works such as \cite{bard2014traveling}, it can be stated as follows:

\begin{align}     \stepcounter{alignnumber}     \setcounter{equation}{0}  \label{MP-1}
    \max\, & \sum_{k\in\mathcal{K}} \sum_{P \in \mathcal{P}_k} \left(Rev_{k,P} - C^T_{k,P} - C^W_{k,P} \right)\lambda_{k,P} \\[3pt]
    \text{s.t.}\, & \sum_{P \in \mathcal{P}_k} \lambda_{k,P} = 1 & \forall k \in \mathcal{K} \\[3pt]
    & \sum_{k\in\mathcal{K}} \sum_{P \in \mathcal{P}_k} a^i_{k,P} \lambda_{k,P} \leq 1 & \forall i \in \mathcal{I}^{\mathcal{R}} \label{MP1con1} \\[3pt]
    & \lambda_{k, P} \in \lbrace 0, 1 \rbrace & \forall k \in \mathcal{K} , P \in \mathcal{P}_k \label{MP1-last}
\end{align}

Objective function~\hyperref[MP-1]{(\ref{MP-1})} maximizes the total revenue. Constraints~\hyperref[MP-1]{(\ref{MP1con1})} ensure that each caregiver is assigned to exactly one plan. Constraints~\hyperref[MP-1]{(\ref{MP1-last})} ensure that new patients are either rejected or accepted and assigned to exactly one caregiver's plan.

\paragraph{Pricing problems.}
The LP dual (\ref{dualobj}) -- (\ref{duallast}) of the master problem's LP relaxation is given in Appendix~\ref{appendix:duals}. Let $v^* \in \mathbb{R}^{|\mathcal{K}|}$ and $w^* \in \mathbb{R}_{+}^{|\mathcal{I}^{\mathcal{R}}|}$ be a solution of the LP dual.  We use the binary parameter $a^i_{R}$ to indicate whether patient $i$ is visited by caregiver $k$ in route $R$.
Then for each caregiver $k \in \mathcal{K}$, we must solve the pricing problem:

\begin{align}  \stepcounter{alignnumber}   \tag{6}   \label{1st-pricing}
    \max_{P \in \mathcal{P}_k} \, & Rev_{k,P} - C^T_{k,P} - C^W_{k,P}   - \sum_{i\in\mathcal{I}^{\mathcal{R}}} w_i^*\, a^i_{k,P}  
\end{align}

The full pricing problem (\ref{pricing}) -- (\ref{pricingend}) is stated in Appendix~\ref{1stlevel_pricing_MIP}
and is seen to be a MILP with $ \mathcal{O}\left(|\mathcal{I}^{\mathcal{R}} \cup \mathcal{I}^{\mathcal{E}}_k|^2\cdot|\mathcal{D}|\right)$ binary variables and $\mathcal{O}\left(|\mathcal{I}^{\mathcal{R}} \cup \mathcal{I}^{\mathcal{E}}_k|\cdot|\mathcal{D}|\cdot |\Gamma^p|\cdot|\Gamma^t|\right)$ continuous variables. However, each pricing problem is further decomposable into days.

\subsubsection{Second-level decomposition.} We further decompose each caregiver's pricing problem into a master problem and $|\mathcal{D}|$ pricing subproblems. The objective of each subproblem is to find a better route for day $d \in \mathcal{D}$. The master problem combines the daily routes into a weekly plan.

\paragraph{Notation.} 
  For any given caregiver, we use the binary variable $\lambda_{d,r}$ to decide whether on day $d$ he or she takes route $r$. The set of feasible routes on day $d \in \mathcal{D}$ is denoted as $\mathcal{R}_d$. The binary parameter $a^i_{d,r}$ indicates whether patient $i$ is treated on day $d$ in route $r$. The revenue generated by a route $r$ executed by caregiver $k$ is denoted by $Rev_{k,r}$, and similarly the travel cost on day $d$ by $C^T_{d,r}$ and wage cost by $C^W_{k,r}$.

\paragraph{Master problem for caregiver $k$.} 

\begin{align}     \stepcounter{alignnumber}     \setcounter{equation}{0}  \label{MP2-first}
    \max\, & \sum_{d\in\mathcal{D}} \sum_{r \in \mathcal{R}_d} \left( Rev_{d, r} - C^T_{d,r} - C^W_{k,r} - \sum_{i\in\mathcal{I}^{\mathcal{R}}} \frac{w_i^*}{v_i}a_{d, r}^i  \right) \lambda_{d,r} \\
    \text{s.t.}\, & \sum_{r \in \mathcal{R}_d} \lambda_{d,r} = 1 & \forall d \in \mathcal{D} \\
    & \sum_{d\in\mathcal{D}} \sum_{r \in \mathcal{R}_d}a^i_{d,r}  \lambda_{d,r} \leq v_i & \forall i \in \mathcal{I}^{\mathcal{R}} \label{MP2con1} \\
     & \sum_{r \in \mathcal{R}_d}v_i a^i_{d,r}   \lambda_{d,r} - \sum_{\mu\in\mathcal{D}} \sum_{r \in \mathcal{R}_\mu}a^i_{\mu,r}  \lambda_{\mu,r}  \leq 0  & \forall d \in \mathcal{D} ,i \in \mathcal{I}^{\mathcal{R}} \label{MP2con2}\\
    & \sum_{r \in \mathcal{R}_d}a_R^i \lambda_{d,r} + \sum_{\mu \in [d_i]_0}\sum_{r \in \mathcal{R}_d}a_R^i\lambda^{d + \mu,r} \leq 1 & \forall d \in \mathcal{D} ,i \in \mathcal{I}^{\mathcal{R}} \label{MP2con3} \\
    & \lambda_{d,r} \in \lbrace 0, 1 \rbrace   & \forall d \in \mathcal{D} ,r \in \mathcal{R}_d \label{MP2-last} 
\end{align}

Objective function~\hyperref[MP-1]{(\ref{MP2-first})} maximizes the total profit and penalizes the acceptance of new patients. 
Constraints~\hyperref[MP-1]{(\ref{MP2con1})} ensure that the caregiver takes exactly one route on each day. Constraints~\hyperref[MP-1]{(\ref{MP2con2})} and~\hyperref[MP-1]{(\ref{MP2con3})} ensure that new patients are either rejected or receive the required number of treatments. Constraints~\hyperref[MP-1]{(\ref{MP2-last})} ensure a minimum number of days between any two visits of a patient.

\paragraph{Pricing subproblems for caregiver $k$ on day $d$.} 
The LP dual (\ref{dual22}) -- (\ref{dual22end}) of the master problem's LP relaxation is given in Appendix~\ref{appendix:duals}. Let $u^* \in \mathbb{R}^{|\mathcal{D}|}$, $z^* \in \mathbb{R}_{\geq 0 }^{|\mathcal{I}^{\mathcal{R}}|}$, $y^* \in \mathbb{R}_{\geq 0 }^{|\mathcal{D}|\times|\mathcal{I}^{\mathcal{R}}|}$, and $q^* \in \mathbb{R}_{\geq 0 }^{|\mathcal{D}|\times|\mathcal{I}^{\mathcal{R}}|}$ be a solution of the LP dual. For better readability, we denote $\mathcal{I}_{k} := \mathcal{I}^{\mathcal{R}} \cup \mathcal{I}^{\mathcal{E}}_k$. Then we must solve the following pricing problem for caregiver $k$ on day $d$:

\begin{align}     \stepcounter{alignnumber}     \setcounter{equation}{0}  \label{MP-pricing}
    \max_{r \in \mathcal{R}_d} \, & Rev_{r} - C^T_{r} - C^W_{r}  \\ 
    &   -  \left( \sum_{i \in \mathcal{I}^{\mathcal{R}}} \frac{w_{i}^*}{v_i} + \sum_{i\in\mathcal{I}^{\mathcal{R}}} z_i^* + \sum_{i \in \mathcal{I}^{\mathcal{R}}} \left( \sum_{\tilde{d} \in \mathcal{D}} y_{i\tilde{d}}^* -  v_i y_{id}^*  \right)
 + \sum_{i\in\mathcal{I}_{k}}\sum_{\substack{\mu\in [|\mathcal{D}| - d_i]: \\ \mu \in [d_i -d]_0}} q_{i,\mu}^*\, \right) a_r^i \label{2nd-pricing}
\end{align}

The corresponding MILP that represents the second-level pricing subproblem is stated in Appendix~\ref{2ndlevel_pricing_MIP}. It has $\mathcal{O}\left(|\mathcal{I}_{k}|^2\right)$ binary variables, $\mathcal{O}\left(|\mathcal{I}_{k}|\cdot |\Gamma^p|\cdot|\Gamma^t|\right)$  continuous variables, and $\mathcal{O}\left(|\mathcal{I}_{k}|^2\cdot |\Gamma^p|\cdot|\Gamma^t|\right)$  constraints.

\subsection{Column generation}

Generating all feasible daily routes $\mathcal{R}_{k,d}$ for (\ref{MP2-first}) -- (\ref{MP2-last})  and all weekly plans $\mathcal{P}_k$ for (\ref{MP-1}) -- (\ref{MP1-last}) explicitly is impractical due to the typically high number of potential routes and plans. Therefore, we employ  column generation to solve the linear relaxation of the master problems \citep{lubbecke2005selected, desrosiers2011branch}.

\subsubsection{Initialization.} \label{Initialization}
Column generation begins with computing a subset $\tilde{\mathcal{P}}_k \subset \mathcal{P}_k$ of weekly plans for each caregiver $k$ and a subset $\tilde{\mathcal{R}}_{k,d} \subset \mathcal{R}_{k,d}$ of routes for each day $d$, respectively. We initialize the subsets by devising two feasible solutions: The first solution entails visiting all existing patients on the pre-determined days while rejecting all new patients. The second solution employs a greedy heuristic to accept new patients.

The greedy heuristic systematically selects promising assignments as follows. An assignment consists of a patient $i \in \mathcal{I}^{\mathcal{R}}$, a compatible caregiver $k$, and a day pattern $D \subset \mathcal{D}$, with $|D| = v_i$ and the required number $d_i$ of days between any two days in $D$. The value of this assignment is the revenue minus an estimation of the accompanying wage  and travel cost: $R_{ki} -  \underaccent{-}{C^W}_{k}\cdot \left( \min_{j \in \mathcal{I}_k \cup \lbrace 0 \rbrace} \bar{t}_{ij} + \bar{p}_i \right) - \min_{j \in \mathcal{I}_k \cup \lbrace 0 \rbrace} C^T_{dji}$. The wage cost for a new patient is approximated by the caregiver's wage cost multiplied by the expected service time, plus the minimum possible expected travel time to reach the patient. The travel cost is estimated as the minimum travel cost from any potential location to the new patient. Given an approximate value for each assignment, we iteratively choose the most promising assignment. We assign the new patient to the corresponding caregiver and day pattern, ensuring the corresponding daily routes remain feasible.

These subsets, $\tilde{\mathcal{P}}_k \subset \mathcal{P}_k$ and $\tilde{\mathcal{R}}_{k,d} \subset \mathcal{R}_{k,d}$, initialized by the reject-all-new-patients solution and the greedy-assignment solution, serve as initial columns for the respective restricted master problems. That is, the first-layer restricted master problem is derived by substituting $\mathcal{P}_{k}$ with $\tilde{\mathcal{P}}_{k}$ in the linear relaxation of the master problem (\ref{MP2-first}) -- (\ref{MP2-last}). Similarly, the second-layer restricted master problem is derived by replacing $\mathcal{R}_{k,d}$ with $\tilde{\mathcal{R}}_{k,d}$ in the linear relaxation of the master problem (\ref{MP-1}) -- (\ref{MP1-last}).

\subsubsection{Pricing problems.} The restricted master problems only include a subset of feasible columns. Consequently, they yield solutions that are feasible for the linear relaxation of the original master problems but not necessarily optimal. We verify optimality through the pricing subproblems.

At the first level, we determine if the optimal solution of the restricted master problem remains optimal for the linear relaxation of the original master problem by determining whether there exists  a weekly plan with a value exceeding the dual objective function value $y^*_k$ for each caregiver $k$, as determined in (\ref{1st-pricing}). This assessment is carried out using the second-level branch-and-price procedure. 

At the second level, we evaluate the optimality of the restricted master problem by checking for each day $d$ whether there exists a route with a value greater than the dual objective function value $u^*_d$ according to (\ref{MP-pricing}) -- (\ref{2nd-pricing}). 
It is worth noting that, owing to revenue maximization, our second-level pricing problems deviate from the common resource-constrained elementary shortest path pricing problems, and align with the more difficult robust longest path problem. To tackle the second-level pricing problems, we employ a simple branch-and-bound procedure, systematically accepting new patients, as outlined in the remainder of this section.

\paragraph{Branch-and-bound procedure.}
The procedure is strongly guided by upper bounds $\bar{v}_{kd}(i)$ on the value generated by caregiver $k$ visiting patient $i \in \mathcal{I}_k$ on day $d$, calculated as:  
\begin{align}     \stepcounter{alignnumber}     \setcounter{equation}{0} 
   \bar{v}_{kd}(i) := R_{ki} & -  \underaccent{-}{C^W}_{k}\cdot \left( \min_{j \in \mathcal{I}_k \cup \lbrace 0 \rbrace} \bar{t}_{ij} + \bar{p}_i \right) - \min_{j \in \mathcal{I}_k \cup \lbrace 0 \rbrace} C^T_{dji} \label{lb1}  \\
    & - \mathds{1}_{i \in \mathcal{I}^{\mathcal{R}}} \left( \frac{w_{i}^*}{v_i} +  z_i^* + \left( \sum_{\mu \in \mathcal{D}} y_{i\mu}^* -  v_i y_{id}^*  \right)
 + \sum_{\substack{\mu\in [|\mathcal{D}| - d_i]: \mu \in [d_i -d]_0}} q_{i\mu}^*\, \right) \label{lb2}
\end{align}
The value $\bar{v}_{kd}(i)$ is derived from the revenue obtained from visiting patient $i$ minus a lower bound on wage and travel costs incurred by by the visit to patient $i$ in (\ref{lb1}), alongside the dual values delineated in (\ref{lb2}). It is noteworthy that the dual values strongly inform the upper bound. 

First, the upper bounds are used for a pre-selection of patients, where new patients with a value $\bar{v}_{kd}(i) \leq 0$ are excluded from consideration in the branch-and-bound procedure. Not visiting them leads to a higher value in the second-level pricing problem. 
Secondly, the upper bounds $\bar{v}_{kd}(i)$ allow for the computation of an upper bound on the objective function value of the second-level pricing problem. Let $n_k$ be an upper bound on the number of patients that caregiver $k$ can visit in a single day, and binary variables $\delta_i$ indicate whether caregiver $k$ visits patient $i \in \mathcal{I}^{\mathcal{R}}$ on day $d$ or not. Then the second-level pricing problem (\ref{MP-pricing}) -- (\ref{2nd-pricing}) is upper bounded by the solution of the following multi-knapsack problem:

\begin{align}     \stepcounter{alignnumber}     \setcounter{equation}{0}  \label{knapsack} 
    \max\, & \sum\limits_{i\in\mathcal{I}^{E|_{k,d}}} \bar{v}_{kd}(i) +    \sum\limits_{i \in \mathcal{I}^{\mathcal{R}}} \bar{v}_{kd}(i)\cdot \delta_i
         \\ \label{constr_workhours}
          \text{s.t.} \, 
              & \sum_{i \in \mathcal{I}^{\mathcal{R}}}  \left( \min_{j \in \mathcal{I}_k \cup \lbrace 0 \rbrace} \bar{t}_{ij} + \bar{p}_i \right)\delta_i \leq \bar{W}_{{ kd} } - \underaccent{\bar}{W}_{{ kd} } - \sum_{i\in\mathcal{I}^{E|_{k,d}}} \left( \min_{j \in \mathcal{I}_k \cup \lbrace 0 \rbrace} \bar{t}_{ij} + \bar{p}_i \right) 
              \\
    & \sum_{i \in \mathcal{I}^{\mathcal{R}}} \delta_i \leq n_{k,d} - |\mathcal{I}^{E|_{k,d}}| \label{constr_maxnumber}\\
   & \delta_i \in \{ 0, 1 \} & \forall i \in \mathcal{I}^{\mathcal{R}} \label{upperboundend}
\end{align}

In this problem, the objective function (\ref{knapsack}) aims to maximize the sum of upper bounds $\bar{v}_{kd}(i)$ for the value generated by visiting new patients. Constraint (\ref{constr_workhours}) ensures that a lower bound on the total workload does not exceed the caregiver's working hours, while constraint (\ref{constr_maxnumber}) ensures that at most $n_k$ patients are visited in total.

Now, if the upper bound provided by (\ref{knapsack}) -- (\ref{upperboundend}) does not surpass the dual value $u^*_d$ or if the problem is infeasible, it provides evidence that no improving column exists. However, if it exceeds the dual value, the assignment determined by the multi-knapsack problem may or may not lead to an improving column. To verify, we solve the corresponding robust traveling salesman subproblems with time windows and wage cost (R-TSPTW) for caregiver $k$ to determine whether the assignment $\mathcal{I}_{TSP} := \mathcal{I}^{E|{k,d}} \cup \left\lbrace i \in \mathcal{I}_k : \delta_i^* = 1 \right\rbrace$ for the current knapsack solution $\delta^*$ indeed yields a better column:

\begin{align}     \stepcounter{alignnumber}     \setcounter{equation}{0}  \label{TSP}
    \min\, & \multicolumn{2}{l}{  $  \sum\limits_{\substack{i, j \in \mathcal{I}_{TSP}}} C^T_{dij}  x_{ij} +   C^W_k    \left(s_{n+1} -\underaccent{\bar}{W}_{{ kd} } \right) $
       }   &\\
     \text{s.t.} \, 
         &  \multicolumn{2}{l}{ $s_{i,\gamma^p, \gamma^t}  \in [ \underaccent{\bar}{T}_{{ id} }, \bar{T}_{{ id} }] \qquad \qquad \qquad \qquad \qquad \qquad \qquad \qquad  \forall i \in \mathcal{I}_{TSP} \, , \gamma^p \in [\Gamma^p]_0 \, , \gamma^t \in [\Gamma^t]_0 \quad \label{TSP_timewindows}$ } &\\
          & \multicolumn{2}{l}{ $s_{i,0,   \zeta^t}   \geq  \left(\underaccent{\bar}{W}_{{ kd} } + \bar{t}_{0i} +   \zeta^t\hat{t}_{0i}\right)  x_{0i}  + (\underaccent{\bar}{W}_{{ kd} } - \bar{T}_{{ id} })(1- x_{i,n+1} )   \qquad\quad \forall i\in\mathcal{I}_{TSP}\,  \zeta^p, \zeta^t \in \lbrace 0, 1 \rbrace$}\qquad \label{TSP_starttime0} & \\
           & 
    \multicolumn{2}{l}{$s_{i,\gamma^p +   \zeta^p, \gamma^t +   \zeta^t}    \geq  s_{j,\gamma^p , \gamma^t}  + (\bar{p}_j +   \zeta^p\hat{p}_j + \bar{t}_{ji} +   \zeta^t\hat{t}_{ji})  x_{ji} $} \nonumber  &\\
    &\multicolumn{2}{l}{ $\qquad\qquad \qquad + ( \underaccent{\bar}{T}_{{ id} } - \bar{T}_{{ jd} } )(1-  x_{ji} ) \qquad\quad\quad\quad\quad \, \forall i,j\in\mathcal{I}_{TSP}  \,, \gamma^p \in [\Gamma^p]_0 \, , \gamma^t \in [\Gamma^t]_0$ } \quad\nonumber &\\ 
    & \multicolumn{2}{l}{ $\qquad\qquad\qquad\qquad\qquad\qquad\qquad\qquad \,\,\,\,\,\, \forall   \zeta^p,   \zeta^t \in \lbrace 0, 1 \rbrace: \gamma^p +   \zeta^p \leq \Gamma^p, \gamma^t +   \zeta^t \leq \Gamma^t \qquad\label{TSP_starttime_i}$ }  &\\
         & \multicolumn{2}{l}{ $s_{n+1}   \geq  s_{i,\Gamma^p - \zeta^p, \Gamma^t - \zeta^t}  + (\bar{p}_i + \zeta^p\hat{p}_i + \bar{t}_{i,n+1} + \zeta^t \hat{t}_{i,n+1})  x_{i,n+1} \qquad \quad\,\,\,\,\,  \forall i\in\mathcal{I}_{TSP}\,,  \zeta^t \in \lbrace 0, 1 \rbrace$} \label{TSP_startime_end} & \\
          & \multicolumn{2}{l}{$s_{n+1}  \in [\underaccent{\bar}{W}_{kd} , \bar{W}_{ kd} ]$  } \label{TSP_workhours} & \\
            &  \multicolumn{2}{l}{$\sum_{j\in\mathcal{I}_{TSP}} x_{ij}  = 1   \qquad \qquad\qquad\qquad\quad \,\,\,\qquad\qquad\qquad\qquad\qquad\qquad\qquad\qquad \forall i \in \mathcal{I}_{TSP} \label{TSP_visitallpatients}    $      } &\\
    & \multicolumn{2}{l}{$\sum_{j\in\mathcal{I}_{TSP}} x_{ji}   -\sum_{j\in\mathcal{I}_{TSP}} x_{ij}  = 0, \sum_{i\in\mathcal{I}_{TSP}} x_{ii}  = 0   \qquad\qquad\qquad\qquad\qquad\qquad \forall i\in \mathcal{I}_{TSP} $}\label{TSP_flowconservation}&  \\
    & \multicolumn{2}{l}{$\sum_{\substack{i \in \mathcal{I}_{TSP}  \cup \lbrace n+1 \rbrace}} x_{0,i}  = 1, \sum_{\substack{i\in \mathcal{I}_{TSP}  \cup \lbrace0 \rbrace}} x_{i,n+1}  = 1$ }   \label{TSP_startenddepot} &  \\
   & \multicolumn{2}{l}{$x_{ij} \in \lbrace 0,1 \rbrace   \,\qquad\qquad\qquad\qquad\qquad\qquad\qquad\qquad\qquad\qquad \forall i, j\in \mathcal{I}_{TSP} \cup\lbrace 0 \rbrace\cup \lbrace n+1 \rbrace$} & \\
   & \multicolumn{2}{l}{$s_{i,\gamma^p, \gamma^t} , s_{n+1} \geq 0   \,\,\,\,\,\,\qquad\qquad\qquad\qquad\qquad\qquad\qquad\qquad \forall i \in \mathcal{I}_{TSP} \,,  \gamma^p \in [\Gamma^p]_0 \, , \gamma^t \in [\Gamma^t]_0 $} \label{lastone}
\end{align}

Objective function~(\ref{TSP}) of this subproblem minimizes the total travel cost plus total wage cost. Constraints~(\ref{TSP_timewindows}) ensure patients are visited within their time windows, even in a worst-case scenario. Constraints~(\ref{TSP_starttime0}) set the first visit's start time with respect to the caregiver's working start time and worst-case travel time. Constraints~(\ref{TSP_starttime_i}) ensure patients are visited after finishing the previous patient's visit and traveling, for all scenarios of service and travel times. Constraints~(\ref{TSP_startime_end}) accordingly define the earliest return time to the depot. Constraint~(\ref{TSP_workhours}) ensures the caregiver returns within working hours. Constraints~(\ref{TSP_visitallpatients}) guarantee that all patients are visited. Constraints~(\ref{TSP_flowconservation}) define flow conservation, and constraints~(\ref{TSP_startenddepot}) specify the caregiver starts and ends at the depot.

The R-TSPTW subproblem has $\mathcal{O}\left( (n_k)^2 \right)$ binary variables, $\mathcal{O}\left( n_k \cdot \Gamma^p \cdot \Gamma^t \right)$ continuous variables, and $\mathcal{O}\left( (n_k)^2 \cdot \Gamma^p \cdot \Gamma^t \right)$ constraints. Efficient solutions can be found using state-of-the-art MIP solvers, as caregivers in home healthcare typically visit at most $n_k = 8$ patients per day, and the uncertainty budgets w.l.o.g. satisfy  $\Gamma^p \leq n_k$ and $\Gamma^t \leq n_k + 2$ (see Section~\ref{maxbudget}). Furthermore, we emphasize that the objective function value does not include any of the dual variables. Hence, we globally store the R-TSPTW routes and their corresponding values, obviating the need to recalculate (\ref{TSP}) -- (\ref{lastone}) for identical patient sets at other nodes of the B\&B search tree. 

Finally, if the R-TSPTW's objective function value exceeds the dual value $u_d^*$, the corresponding route yields an improving column; otherwise, we branch on the most promising patient with respect to the upper bound $\bar{v}_{kd}(i)$ in order to either find an improving assignment or prove that none exists. Thus, focusing on promising routes, our branch-and-bound procedure efficiently explores the solution space of daily routes.

\subsection{Branching}

At each level, the column generation procedure yields an optimal solution for the linear relaxation of the corresponding master problem, which may not be integer. To obtain an optimal integer solution, we integrate the column generation procedure at each level with a branch-and-bound tree, resulting in a nested branch-and-price method.

\begin{algorithm}[b]
\caption{Overview of first-level branching}\label{first_alg}
\For{$i \in \mathcal{I}^{\mathcal{R}}$ not yet branched upon}{
    $\rm{UB}_{i, \emptyset} \leftarrow 0 $\;
    \For{$ k \in \mathcal{K}$}{ 
    $\rm{UB}_{i,  k} \leftarrow \max_{D \in \rm{day\textunderscore patterns}_{ik}}\rm{Val}(D) $\;}
    $i\textunderscore \rm{branch} \leftarrow \argmax_i \left\lbrace \max_{k} \left\lbrace \rm{UB}_{i, k} \right\rbrace \right\rbrace$\;
}
\KwRet{$i\textunderscore \rm{branch}$, $\rm{caregivers}_{i\textunderscore \rm{branch},k}$ 
}
\end{algorithm}

Upon termination of the column generation procedure with a solution $\lambda$, we check if the solution already constitutes an integer solution. If yes, we prune the node with this integer solution. Otherwise, we apply the following branching strategy: In the first level, patients are branched upon caregivers. In the second level, patients are branched upon day patterns. At each level, we select the most promising remaining patient. 
\begin{algorithm}[t]
\caption{Overview of second-level branching}\label{second_alg}
\For{$i \in \mathcal{I}^{\mathcal{R}}_k$ not yet branched upon}{
    $\rm{UB}_{i, \emptyset} \leftarrow 0 $\;
    \For{$ D \in \rm{day\textunderscore patterns}_{ik}$}{ 
    $\rm{UB}_{i, D} \leftarrow \rm{Val}(D) $\;}
    $i\textunderscore \rm{branch} \leftarrow \argmax_i \left\lbrace \max_{D} \left\lbrace \rm{UB}_{i, D} \right\rbrace \right\rbrace$\;
}
\KwRet{$i\textunderscore \rm{branch}$, $\rm{day\textunderscore patterns}_{i\textunderscore \rm{branch}}$ 
}
\end{algorithm}

Specifically, for each patient $i$ and each caregiver $k$, the set of feasible day patterns, denoted as $\rm{day}\textunderscore \rm{patterns}_{ik}$, encompasses all subsets of days satisfying the correct number of visits per patient and the minimum number of days between any two visits. Optionally, it may ensure that visiting patient $i$ on all days in the pattern does not exceed caregiver $k$'s workload. The value of a day pattern is estimated using the upper bound $Val(D) := \sum_{d \in  D} \left( R_{ki} - C^W_{k,d}\left( \bar{p}_i + \frac{1}{2}\left(\min_{j' \neq j'' \in \mathcal{I}^{E}_{k,d}}{\bar{t}_{j',i}} + {\bar{t}_{i,j''}}\right) \right) - \frac{1}{2}\left(\min_{j'\neq j'' \in \mathcal{I}^{E}_{k,d}}{C^T_{d,j',i}} + {C^T_{d,i,j''}}\right) \right)$. The patient selected for branching is chosen to maximize the estimated value $\rm{UB}$ of the fixed assignment, as described in Procedure~\ref{first_alg} and Procedure~\ref{second_alg}. Furthermore, the corresponding assignment of maximum value is selected for node selection.

\section{Computational Study}
\label{computational_study}

In this section, we evaluate the results obtained with our nested branch-and-price procedure using real-world data obtained from a prominent US healthcare agency based in Dallas, Texas [see \cite{guo2023three})]. We address the agency's task of computing routes and schedules for the forthcoming week. At the end of each day of the current week, the agency decides whether to accept the patient requests received that day. The acceptance decisions are made simultaneously with assignment, routing, and scheduling decisions, while ensuring continuity of care for existing patients. We review the resulting weekly plan and derive management insights about different degrees of uncertainty and time windows.

\subsection{Instances}
In all, we solved 144 test instances generated from the agency's dataset (see Appendix~\ref{data_description} for details).  These instances encompass 5,999 routine visits conducted by 36 caregivers during a four-week period from July 23 to August 20, 2018. The first week is set apart to establish an initial weekly plan using the greedy heuristic described in Section~\ref{Initialization}, while the subsequent three weeks are used to compute and evaluate optimal weekly home healthcare plans with the nested B\&P procedure.

Section~\ref{fixed} outlines the fixed parameter values that remain constant across all instances. Section~\ref{varied} introduces the factors—discipline, regional branch, time window size, and uncertainty budget—that we systematically vary using a full factorial design. 

 \subsubsection{Fixed factors.} \label{fixed} The parameters—revenue and cost, along with availability of caregivers and new patients—remain constant across all instances.

\paragraph{Revenue \& cost.}
The revenue per visit is calculated using the nursing hourly billing rate of 102 US\$, as reported by the home healthcare agency under consideration. We assume that the services provided by  registered nurses, occupational therapists, and speech therapists generate exactly this hourly revenue, while a 10\% discount is applied to licensed practical nurses, skilled nurses, occupational therapist aides, and physical therapy assistants, and a 10\% surcharge is applied for physical therapy. For hourly-paid caregivers, we refer to Table~\ref{tab:services} for wage cost. If caregivers are paid per visit, we subtract the wage cost, proportional to the expected service time, from the revenue of each visit. Salaried caregivers incur a fixed wage cost in this operational decision, which is the weekly number of work hours multiplied by the hourly wage cost. Travel cost is estimated at 0.56 US\$ per mile, based on the 2021 IRS standard mileage rate.
 
\paragraph{Caregivers.}
We establish a fixed group of caregivers across all weeks, determined based on the week of August 1, 2018. Caregivers are available from Monday to Friday, with full-time working hours from 9 am to 5 pm and part-time working hours from 8 am to 12 pm. The employment type and compensation details are presented in Table~\ref{tab:caregivers} in Appendix~\ref{results_description}. All caregivers with a contract or working on an as-needed basis (PRN) are assumed to be available full-time. Each caregiver is assumed to be proficient in providing any service that represented at least 5\% of the visits delivered by that caregiver in the historical data.

\paragraph{Patients.}
Patients are classified as either existing or new, depending on whether they were visited in the previous week. Existing patients retain their visits by the same caregiver on the same days. For new patients, however, the assignment needs to be determined. New patients are considered available every day of the week and compatible with any caregiver capable of providing the required service. The number of weekly visits required for each new patient is determined by their mean weekly visits in the historical data. The minimum number of days between routine visits is based on the minimum observed days in the data for each service. For all services this is 0 (but at most 1 visit per day), except for registered nursing, where it is 1 day. 

 \begin{table}[b]
    \centering \footnotesize
    \caption{Services in the data set}
    \begin{tabular}{lrrrrrrrrrrrrrrrrrrrrrrrrr} \toprule
         \bf Service    && \bf Mean frequency &&  \bf Mean treatment period && \bf Expected && \bf {Wage}     \\    
         && (per patient and week) &&   (per patient) && \bf duration && (hourly)   \\  
         \midrule
             \bf SN  &&  1 visit(s)  && 5 weeks &&   50 min && \$26.86    \\
            \bf RN  && 1 visit(s) &&  2 weeks &&  55 min  && \$42.80    \\ 
            \bf LPN  &&  1 visit(s) && 4 weeks &&  50 min &&  \$26.86  
       \\  \midrule
        \bf PTA  && 2 visit(s) &&  5 weeks &&  50 min   && \$31.01      \\ 
           \bf PT   &&  1 visit(s) &&  3 weeks&&  55 min && \$47.10  \\ 
         \bottomrule
    \end{tabular}
    \label{tab:services1}
\end{table}

\subsubsection{Varied factors.} \label{varied} In order to evaluate the robustness of our home healthcare plans, we introduce uncertainty to each instance and then systematically vary the following factors in a full factorial design: discipline, region, time window size, and uncertainty level. In the following, we provide a detailed description of each of these factors.

\paragraph{Disciplines.}
We explore two disciplines: nursing and physical therapy. Each discipline consists of specific services: Nursing includes skilled nursing (SN), registered nursing (RN), and licensed practical nursing (LPN) services. Physical therapy comprises physical therapy (PT) and physical therapy assistance (PTA). By default, a test instance represents one discipline, encompassing all routine visits within that discipline. 
The disciplines and services vary in terms of mean weekly visits per patient, average visit duration, and wage cost. An overview of these variations is provided in Table~\ref{tab:services1}.

\begin{figure}[t]
    \centering
    \includegraphics[scale = 0.85]{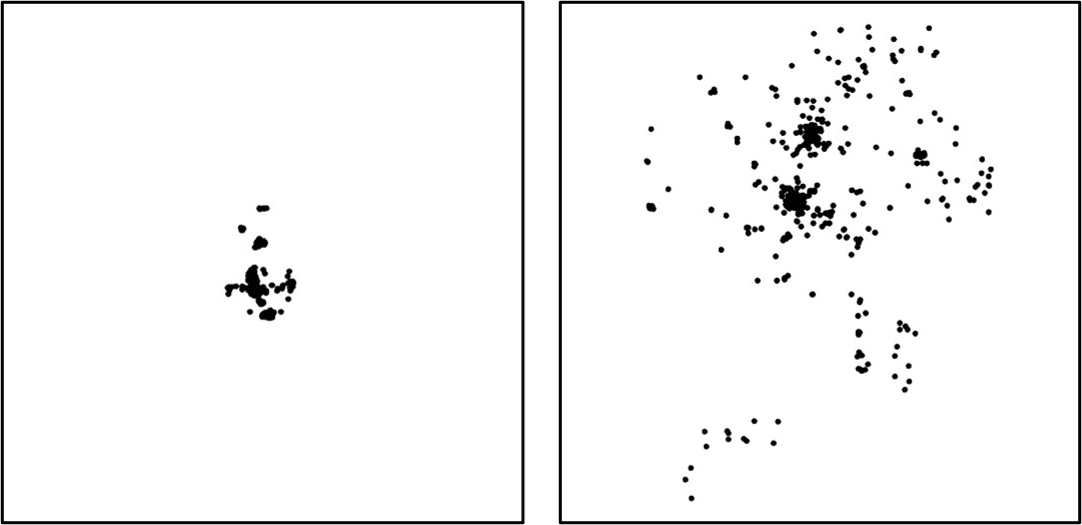}
    \caption{Scatter plots of patient locations in Region-T20 and Region-T60 \qquad \qquad \qquad \qquad \qquad \qquad \qquad \qquad \qquad \qquad \qquad \qquad \qquad \qquad \qquad \qquad \qquad \qquad \qquad \qquad \qquad \qquad \qquad \qquad   \textmd{(arranged in increasing order with respect to average travel time, identical scaling)}}
    \label{fig:regions1}
\end{figure}

\paragraph{Regions.} The two regional branches considered in our computational study differ in terms of population, patient count, average travel time and distance between any two patient locations. We named the branches based on the average travel time, computed using the OSRM routing engine \citep{luxen-vetter-2011}. Specifically, Region-T20 has an average travel time of 20 minutes (and 22.60 km) 
and Region-T60 has an average travel time of 60 minutes (and 74.69 km) between any two patient locations. A scatter plot depicting the patient locations in each region is shown in Figure~\ref{fig:regions1}. Region-T20 and Region-T60 are rural areas with a population below 100,000. 

\paragraph{Time windows.} For existing patients, we explore three sizes of time windows: an exact appointment time (size: 0 minutes), a one-hour time window (size: 60 minutes), and a AM/PM/None time window categories based on the data set.

\paragraph{Uncertainty.} We investigate four levels of uncertainty ranging from no uncertainty to anticipated deviations for all time estimates.
We consider that patients and paths typically necessitate up to a 20\% deviation from their expected service or travel time, respectively. We vary the uncertainty budgets, planning robustly against deviations of none, a quarter, half or all of the expected service and travel times.  
Table~\ref{tab:ceterisparibus} summarizes these levels, characterized by the uncertainty budget $\Gamma^p$ for service times, maximum service time deviations $\hat{p}$ (as a fraction of the expected service times $\bar{p}$), and the uncertainty budget $\Gamma^t$ for travel times, alongside maximum travel time deviations $\hat{t}$ (as a fraction of the expected travel times $\bar{t}$).

\begin{table}[h]
    \centering \footnotesize
    \caption{Levels of uncertainty}
    \begin{tabular}{cccccccccccccccccccccccccc} \toprule
        \multicolumn{2}{c}{None}   && \multicolumn{2}{c}{Low}  && \multicolumn{2}{c}{Medium}  && \multicolumn{2}{c}{High} \\ 
         \cmidrule{1-2}  \cmidrule{4-5} \cmidrule{7-8}   \cmidrule{10-11}   
        $\Gamma^p = 0$, & $\Gamma^t = 0$  && $\Gamma^p = 2$, & $\hat{p} = 20\%$  && $\Gamma^p = 4$, & $\hat{p} = 20\%$ && $\Gamma^p = 8$, & $\hat{p} = 20\%$  \\
         $\hat{p} = 20\%$, & $\hat{t} = 20\%$  && $\Gamma^t = 2$, & $\hat{t} = 20\%$  && $\Gamma^t = 4$, & $\hat{t} = 20\%$ && $\Gamma^t = 8$, & $\hat{t} = 20\%$   \\
 \bottomrule
    \end{tabular}
    \label{tab:ceterisparibus}
\end{table}

\subsection{Runtimes}
All codes were written in Python and all computations performed in a Microsoft Windows 10 environment, running on an Intel Core i7-10510U 1.8 GHz processor with 16 GB memory operating in 64-bit mode. The second-level pricing problems were solved with Gurobi version 9.5.0. The implementation adheres to the specifications outlined in Section~\ref{method}. 

Table~\ref{tab:results_runtime} in Appendix~\ref{results_description} displays the runtimes for the varied factors; i.e., region, discipline, uncertainty level, and time window size. The average runtime is 591.01 seconds. Based on a correlation analysis, we can conclude that the following parameters primarily contribute to the runtime disparities.

\paragraph{Number of caregivers and patients.} The factors impacting runtime most significantly are patient and caregiver counts, along with the daily number of new patient requests. This corresponds to the number of possible decisions, which increases exponentially with each of these factors. Table~\ref{results_runtimes_n} summarizes the runtime effects. 

\begin{table}[b] 
\centering \footnotesize
\caption{Runtimes by the total number $n$ of considered patients, of which $n_{new}$ patients are new \normalfont{(in seconds)}}
\label{results_runtimes_n}
\begin{tabular}{llrrrrrrrrrrrrrrr}
\toprule
 && $n_{new} = 15$ (3 per day) && $n_{new} = 25$ (5 per day)   \\ 
\cmidrule{5-5}  \cmidrule{3-3}
n = 1 - 24: && 10.94
 && 38.87 && (K = 3.00)  \\
n = 25 - 49: && 1.23 && 179.34 &&  (K = 4.90) \\
n = 50 - 74: && 5.09 && 46.65 &&  (K = 6.19) \\
n = 75 - 99: && 75.70 && 158.62 && (K = 6.42) \\
n = 100 - 124: && 201.51  && 346.89 && (K = 6.08) \\
n = 125 - 149: && - && 2\,922.67  && (K = 6.00) \\
\cmidrule{5-5}  \cmidrule{3-3}
  Average:  && \textcolor{green!33!black}{tba} && 380.92  \\
\bottomrule
\end{tabular}
\end{table}

\paragraph{Uncertainty level.} 
Higher robustness requires significantly longer runtimes, as shown in Table~\ref{results_runtimes_nlevel}. Two main observations emerge:
Firstly, there is a significant increase in runtimes up to a medium level of uncertainty, followed by a decrease. This pattern arises because setting the uncertainty budget equivalent to the number of daily patients (or paths) per caregiver mirrors the scenario of no uncertainty, where all data necessitates the maximum deviation (refer to Section~\ref{maxbudget}). Therefore, runtimes peak when the uncertainty budget matches half the daily patients (and paths) per caregiver. 
Secondly, within the nested B\&P procedure, uncertainty primarily impacts the second-level R-TSPTW pricing problem (see (\ref{MP-pricing}) -- (\ref{2nd-pricing})). Since we used Gurobi to solve it, exploring faster solution methods for the R-TSPTW in future research appears promising to further improve runtimes under uncertainty \citep{bartolini2021robust, lera2022dynamic}.

\begin{table}[t] 
\centering \footnotesize
\caption{Runtimes by the number $n$ of patients and the level of robustness \normalfont{(in seconds)}}
\label{results_runtimes_nlevel}
\begin{tabular}{llrrrrrrrrrrrrrrr}
\toprule
 && None && Low && Medium && High  \\ 
\cmidrule{5-5} \cmidrule{7-7} \cmidrule{3-3} \cmidrule{9-9} 
n = 1 - 24: && 8.92 &  & 61.25 &  &145.43 && 5.27  \\ 
n = 25 - 49: &&  50.72 &  & 255.42 &  & 708.54 && 15.42  \\ 
n = 50 - 74: && 2.08 &  & 1.71 &  &  87.03 && 5.22  \\ 
n = 75 - 99: && 5.02 &  & 15.81 &  & 254.51 && 549.68  \\ 
n = 100 - 124: && 37.65 &  & 305.87 &  &  1\,017.03 && 148.43  \\ 
n = 125 - 149: && 170.59 &  & 459.46 &  & 10\,067.92 && 305.24  \\
\cmidrule{5-5} \cmidrule{7-7} \cmidrule{3-3} \cmidrule{9-9} 
  Average:  && 45.83 &  & 183.25 &  &  2\,046.74 && 88.21   \\
\bottomrule
\end{tabular}
\end{table}

\paragraph{Time windows.} Time window size also impacts runtimes, albeit to a smaller extent compared to patient and caregiver counts. Tight time windows for existing patients (0 minutes) result in an average runtime of 462.94 seconds, narrow time windows (60 minutes) 42.27 seconds, and wide time windows (AM/PM/None) result in a runtime of 34.17 seconds. This decrease in runtimes with the time window size is primarily due to the number of feasible routes, which decreases significantly with shorter time windows.

\paragraph{Other.} 
Finally, we identify factors with relevant but less significant effects on runtime based on the average correlation of all parameters when fixing the varied factors. As the week progresses and more capacity is used, the procedure accelerates (runtime correlation is -0.13 with the day index). A slight decrease in runtime occurs with a higher number $v$ of weekly visits per patient, which we explain by the number $ \binom{5}{v} $ of day patterns (correlation of -0.13). We also found that the optimality gap between the optimal solution and the first heuristic solution to have a negligible effect on runtime. Not surprisingly, these observations collectively suggest that faster runtimes are associated with fewer feasible solutions.

\subsection{Managerial insights}
We investigate the impact of robustness and time window size on optimal weekly home healthcare plans, focusing on managerial insights. In Section~\ref{profsens}, we analyze profit sensitivity, while Section~\ref{decisions} explores optimal acceptance, assignment, routing, and scheduling decisions. For a comprehensive overview of results across various regions, disciplines, time windows, and uncertainty levels, we refer to Table~\ref{tab:results_summary} in Appendix~\ref{results_description}. 

\begin{table}[b]
\centering 
\footnotesize
\caption{Profit by level of robustness and time window size \textmd{(in USD)}}
\begin{tabular}{llrrrrrrrrrrrrrrr}
\toprule
&& None && Low && Medium && High \\
\cmidrule{5-5} \cmidrule{7-7} \cmidrule{3-3} \cmidrule{9-9} 
Wide: && 5,624 && 5,361  && 5,095 && 4,880 \\
Narrow: && 5,297 && 5,026  && 4,825 && 4,767 \\
Tight: && 5,297 && 4,796 &&  4,735 && 3,306 \\
\cmidrule{5-5} \cmidrule{7-7} \cmidrule{3-3} \cmidrule{9-9} 
Average: && 5,406 && 5,061 && 4,886 && 4,318 \\
\bottomrule
\end{tabular}
\label{tab:profit}
\end{table}

\subsubsection{Profit sensitivity.} \label{profsens}
 Across all instances, the weekly profit averages \$4,964.46, predominantly driven by a revenue of \$6,805.32, surpassing average weekly wage and travel costs of \$1,438.68 and \$402.18, respectively. The profitability of each visit, averaging \$96.97, emerges as the primary determinant of profit. Our subsequent analysis examines the effects of robustness and time window size on profit, revenue, wage cost, and travel cost, as outlined in Table~\ref{tab:profit} to Table~\ref{tab:travelcost}.
\begin{wrapfigure}{r}{0.33\textwidth}
    \centering \small 
    \begin{subfigure}[b]{0.33\textwidth}
        \centering \small
        \begin{tikzpicture}[scale=1, transform shape]
            \node[align=center] (legend) at (2,6.5) {
            \ref{pgfplots:plot1} profit, 
            \ref{pgfplots:rev} revenue \\
            \ref{pgfplots:wage} wage cost, 
            \ref{pgfplots:travel} travel cost
        };
            \pgfplotsset{set layers}
            \begin{axis}[
                height=5cm,
                width = 4cm,
                scale only axis,
                x axis line style = {opacity=0},
		axis x line* = bottom,
		axis y line = left,
		enlarge y limits,
		ymajorgrids,
  ylabel={USD},
                xmin=0,xmax=5,
                xtick={1,2,3,4},
                xticklabels={None, Low, Medium, High},
                xtick style = {draw=none},
                xticklabel style={rotate=90, anchor=north east},
                axis y line*=left,
                ylabel style={align=center},
            ]
            \footnotesize
            \addplot[black, line width=0.5pt, mark=*, mark size=0.66pt] coordinates {
                (1, 5406.535730)
                (2, 5061.253212)
                (3, 4885.475768)
                (4, 4353.009758)
            };
                    \label{pgfplots:plot1}     
            \addplot[green!66!black, dashed, line width=0.33pt, mark=*, mark size=0.66pt] coordinates {
                (1, 7264.016334)
                (2, 6896.916464)
                (3, 6714.567816)
                (4, 6232.551561)
            };
             \label{pgfplots:rev}  
            \addplot[TUMorange!66!red, dashed, line width=0.33pt, mark=*, mark size=0.66pt] coordinates {
                (1, 1438.680000)
                (2, 1438.680000)
                (3, 1438.680000)
                (4, 1488.289655)
            };
            \label{pgfplots:wage}
            \addplot[brown, dashed, line width=0.33pt, mark=*, mark size=0.66pt] coordinates {
                (1, 418.8006038)
                (2, 396.9832526)
                (3, 390.4120473)
                (4, 391.252148)
            };
            \label{pgfplots:travel}
            \end{axis}
        \end{tikzpicture}
        \label{fig:plot2}
    \end{subfigure}
    \caption{\small \bf \textup  Levels of robustness \\}
        \hfill
    \label{fig:three_plots}
\end{wrapfigure}
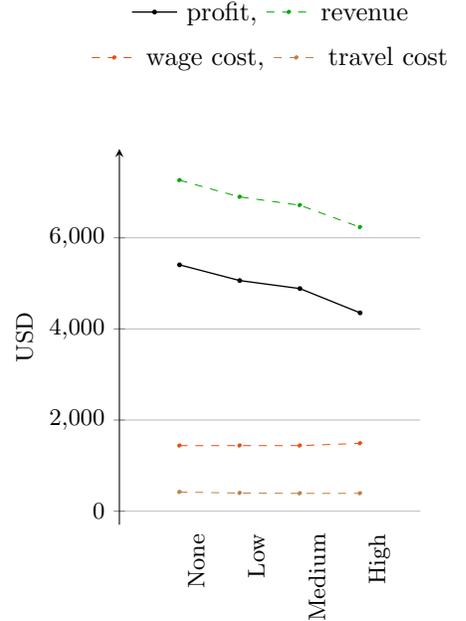

\paragraph{Robustness.} Figure~\ref{fig:three_plots}  illustrates the profit, revenue, and costs across the different levels of robustness. We draw two key observations:
\begin{table}[t]
\centering \footnotesize
\caption{Revenue by level of robustness and time window size \textmd{(in USD)}}
\label{tab:rev}
\begin{tabular}{llrrrrrrrrrrrrrrr}
\toprule
&& None && Low && Medium && High \\
\cmidrule{5-5} \cmidrule{7-7} \cmidrule{3-3} \cmidrule{9-9}
Wide: && 7,496 && 7,201 &&  6,916 && 6,709 \\
Narrow: && 7,151 && 6,860 && 6,661 && 6,598 \\
Tight: && 7,145 && 6,630 && 6,566 && 5,297 \\
\cmidrule{5-5} \cmidrule{7-7} \cmidrule{3-3} \cmidrule{9-9} 
Average: && 7,264 && 6,897 &&  6,715 && 6,201 \\
\bottomrule
\end{tabular}
\end{table}

Firstly, in this case study, the almost constant disparity between revenue and profit curves indicates that robustness primarily decreases profit due to reduced patient acceptance. This is also evident from the observation that, on average, weekly plans with no uncertainty involve caregivers visiting on average $4.34$ patients per day, while low, medium, and high uncertainty permit $4.01$, $3.84$, and $3.75$ patients per day, respectively. 
This  underscores the obvious yet crucial insight that fewer patient visits per day improve operational smoothness. 
However, it is also important to note that robustness can compromise health access if caregiver numbers cannot be augmented concurrently.

Secondly, the transition from no to low robustness decreases revenue by 5.05\%, to medium robustness by 7.56\%, and to high robustness by 14.20\%. In particular, the decline in revenue is small between low and medium robustness. This finding suggests that in some cases, weekly plans can be robustified by adjusting schedules without compromising revenue, see Section~\ref{decisions}.

\begin{table}[t]
\centering \footnotesize
\caption{Wage cost by level of robustness and time window size \textmd{(in USD)}}
\begin{tabular}{llrrrrrrrrrrrrrrr}
\toprule
&& None && Low && Medium && High  \\
\cmidrule{5-5} \cmidrule{7-7} \cmidrule{3-3} \cmidrule{9-9} 
Wide: && 1,439 && 1,439 &&  1,439 && 1,439 \\
Narrow: && 1,439 && 1,439 &&  1,439 && 1,439 \\
Tight: && 1,439 && 1,439 &&  1,439 && 1,599 \\
\cmidrule{5-5} \cmidrule{7-7} \cmidrule{3-3} \cmidrule{9-9} 
Average: && 1,439 && 1,439 &&  1,439 && 1,492 \\
\bottomrule
\end{tabular}
\label{tab:wagecost}
\end{table}

\begin{table}[b]
\centering \footnotesize
\caption{Travel cost by level of robustness and time window size \textmd{(in USD)}}
\begin{tabular}{llrrrrrrrrrrrrrrr}
\toprule
&& None && Low && Medium && High \\
\cmidrule{5-5} \cmidrule{7-7} \cmidrule{3-3} \cmidrule{9-9}
Wide: && 433 && 401 && 382 && 390 \\
Narrow: && 415 && 395  && 397 && 392 \\
Tight: && 409 && 395  && 392 && 392 \\
\cmidrule{5-5} \cmidrule{7-7} \cmidrule{3-3} \cmidrule{9-9} 
Average: && 419 && 397  && 390 && 391 \\
\bottomrule
\end{tabular}
\label{tab:travelcost}
\end{table}

\paragraph{Time windows.} On average, across all instances, wide time windows result in a profit of \$5,412.46, narrow time windows in \$4,922.69, and tight time windows in \$4,540.80. However, the impact of time windows is intricately linked with the level of robustness and warrants a detailed examination.

With no robustness, only a single scenario is considered, hence a single service time is determined. 
Wide time windows yield a profit of \$5,624.34, comprising a revenue of \$7,495.96, with wage and travel costs at \$1,438.68 and \$432.94, respectively.  Narrow time windows induce more rejections, leading to a 4.60\% decrease in revenue, and a 5.80\% profit decline, respectively. Tight time windows only exacerbate this effect slightly, resulting in a 5.82\% decrease in profit compared to wide time windows.

Introducing robustness leads to scenario-dependent service start times. Both the earliest and latest possible start times must fall within a patient's time window. Now, at a medium level of robustness, wide time windows yield a profit of \$5,095.61, a revenue of \$6,916.48, with wage and travel costs at \$1,438.68 and \$382.18, respectively. Once again, narrow time windows lead to more rejections, causing a comparable revenue decrease of 3.69\%, and a profit decrease of 5.30\%. However, under uncertainty, tight time windows exacerbate this trend further, resulting in a revenue decline of 5.07\%, and a profit loss of 7.08\% compared to wide time windows.

Overall, transitioning from wide time windows and no robustness to narrow time windows coupled with a medium level of robustness results in an average profit decrease of 14.20\% in this case study. Primarily, this decline can be attributed to a reduced number of visits. However, the impact on profit fluctuates between different levels of robustness and time windows, suggesting a non-linear price of robustness and time window size.

\subsubsection{Optimal decisions.} \label{decisions}
We now investigate optimal decisions for acceptance, assignment, routing, and scheduling. In contrast to many other approaches in the literature, our methodology makes these decisions simultaneously.

\begin{figure}[htpb]
\centering \footnotesize
\label{tab:rev}
\resizebox{\textwidth}{!}{
\begin{tabular}{p{2.25cm}lrrrrrrrrrrrrrrr}
 && No && Low && Medium && High  \\ 
  && robustness && robustness && robustness && robustness  \\ 
\cmidrule{5-5} \cmidrule{7-7} \cmidrule{3-3} \cmidrule{9-9} 
Wide \qquad \quad time windows:  && \includegraphics[width=5cm]{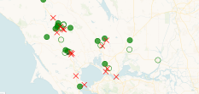} && \includegraphics[width=5cm]{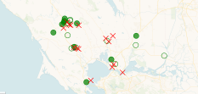} && \includegraphics[width=5cm]{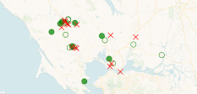} && \includegraphics[width=5cm]{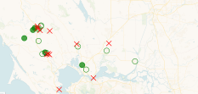}   \\
Narrow \qquad \quad time windows: && \includegraphics[width=5cm]{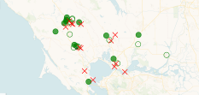} && \includegraphics[width=5cm]{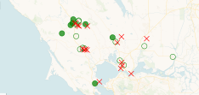} && \includegraphics[width=5cm]{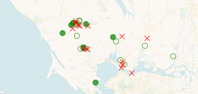} && \includegraphics[width=5cm]{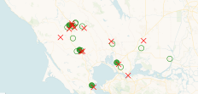}  \\
Tight \qquad \quad time windows: && \includegraphics[width=5cm]{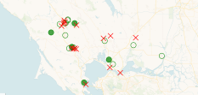} && \includegraphics[width=5cm]{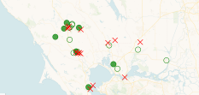} && \includegraphics[width=5cm]{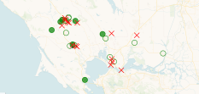} && \includegraphics[width=5cm]{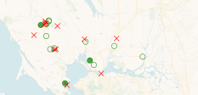} \\
\end{tabular}
}
\caption{Acceptance decision by robustness level and time window size \textmd{(for Region T-60, Nursing)} \\}
 \label{results_accept_levels}
\end{figure}

\paragraph{Acceptance.} We first investigate the number of accepted patients as well as potential selection criteria across varying levels of robustness and time window sizes. Figure~\ref{results_accept_levels} illustrates patient acceptance and rejection for nursing in Region T-60.

If capacity permits, all patients are accepted, as the revenue from each patient typically exceeds the associated wage and travel costs. Selection only occurs when capacity constraints necessitate rejecting patients. On average, 1.53 out of 5 new patients are accepted on day 1, 1.67 on day 2, 1.24 on day 3, 0.89 on day 4, and 0.60 on day 5, indicating a high capacity utilization for the home healthcare agency under consideration. As already indicated in terms of revenue changes, acceptance rates vary with different levels of robustness and time window sizes. For wide time windows, 6.27 patients are typically accepted, compared to 6.04 with narrow time windows, and 5.67 with tight time windows. Similarly, on average, 7.65, 6.17, 5.64, and 4.92 out of 25 new patients are accepted per week with no, low, medium, and high budget robustness, respectively. This represents a decrease in acceptance rates of 19.35\%, 26.27\%, and 35.69\% for low, medium, and high robustness levels compared to no robustness.

Furthermore, the simple greedy heuristic from Section~\ref{Initialization} demonstrates strong performance, achieving an average gap of 7.73\% compared to the optimal solution, relative to the difference between the optimal solution value and the solution value when rejecting all patients, and averaged over all days where the optimal solution accepts at least one patient. Additionally, Figure~\ref{results_accept_levels} illustrates that new patients may be selected even if they live farther away than other candidates, indicating that geographic proximity is not necessarily critical. These observations suggest that a patient's contribution margin may take priority over geographical proximity in the decision-making process. Accordingly, we recommend selecting patients similar to the greedy heuristic, iteratively based on their potential value, determined by revenue minus an estimate of travel and wage costs.

\paragraph{Assignment.} The decision of assigning new patients to caregivers and days carefully balances  maximizing total revenue and minimizing total travel time and distance, while ensuring feasibility with respect to all patient time windows and work hours, even under uncertainty.

\begin{figure}[htbp] 
\footnotesize
\label{tab:rev}
    \begin{tabular}{c} 
\includegraphics[width=4cm, trim=4cm 0.1cm 4cm 1cm, clip]{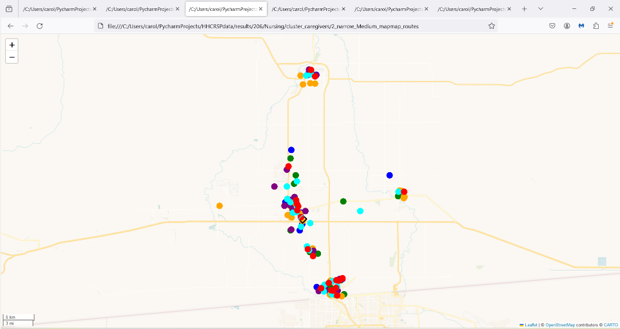} \\
Week   \\ \\
\end{tabular}
\vfill

\begin{subfigure}{\textwidth}
\centering
\resizebox{\textwidth}{!}{
    \begin{tabular}{ccccccccccccccccccccccccc}
\toprule 
\\
\includegraphics[width=4cm,trim=4cm 0.1cm 4cm 1cm, clip] {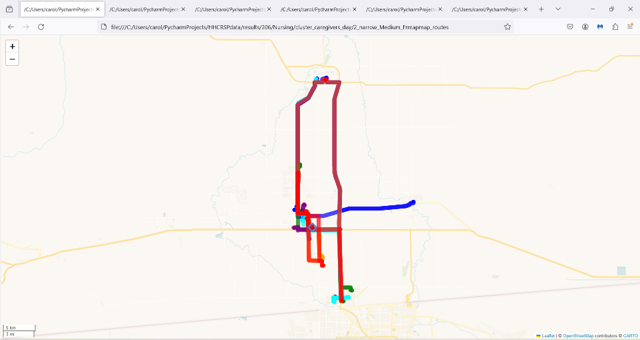} && \includegraphics[width=4cm,trim=4cm 0.1cm 4cm 1cm, clip]{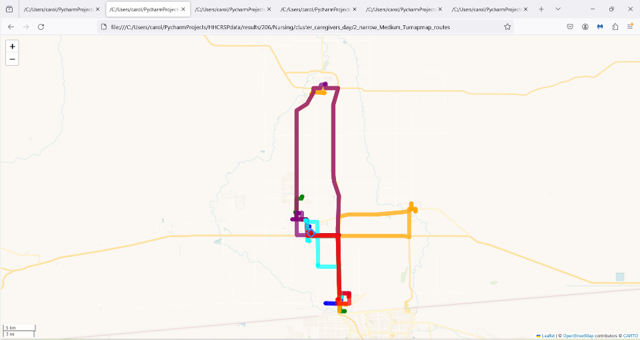}  && \includegraphics[width=4cm,trim=4cm 0.1cm 4cm 1cm, clip]{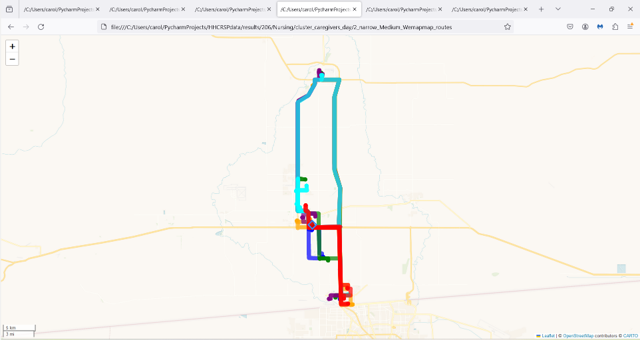} && \includegraphics[width=4cm,trim=4cm 0.1cm 4cm 1cm, clip]{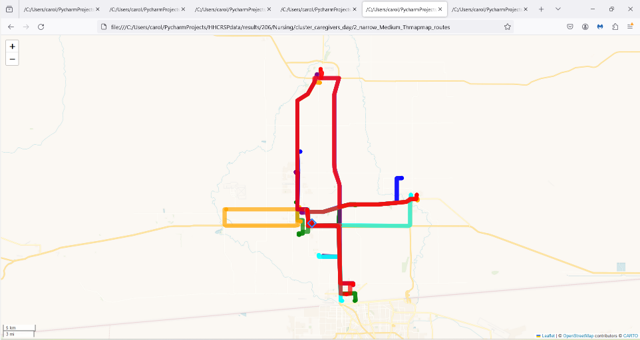} && \includegraphics[width=4cm,trim=4cm 0.1cm 4cm 1cm, clip]{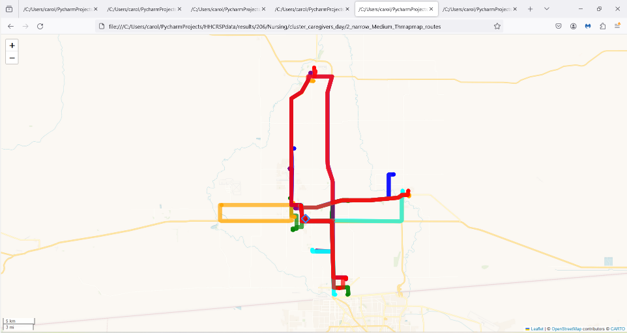} \\
 Day 1 && Day 2 && Day 3 && Day 4 && Day 5  \\ 
 \\
\end{tabular}
}
\end{subfigure}
\caption{Illustration of assignment of patients to caregivers by week and by day \textmd{(for Region T-20, Nursing)} }
\label{results_clusters}
\end{figure}

We find that geographic proximity holds significance primarily in settings with wide time windows. However, even in such cases, there is limited support for the traditional approach of pre-clustering caregivers by region. Instead, we observe geographic clusters by caregiver for each day independently, as depicted in Figure~\ref{results_clusters}. To clarify this observation, let us consider two extreme cases theoretically: one where all patients require one visit per week and another where some patients require visits on all days. In the former case, caregivers visit similar regions daily, aligning with pre-clustering regions by caregiver. Conversely, in the latter case, where patients need only one visit per week, freely assigning a new region for each caregiver each day may be more profitable. Thus, the conventional approach of pre-clustering caregivers by region seems advantageous only when time windows are wide and patients require multiple visits per week.

\vfill

In settings with narrow or tight time windows, we cannot observe distinct patterns based on geographic proximity. We explain this with feasibility considerations. Even if nearby patients cannot be visited within their time windows, it remains more profitable to visit distant patients than none at all. As a result, geographical factors appear less influential with shorter time windows, further compounded by higher levels of robustness. In such situations, we discourage the use of clustering heuristics and suggest either computing optimal solutions or resorting to alternative heuristics like a simple greedy heuristic that iteratively fixes the most promising patient-caregiver-day patterns.

\paragraph{Routes and schedules.} 
Finally, we assess the daily routes and schedules across various levels of robustness and time window sizes. If fewer visits can be accommodated due to higher robustness or shorter time windows, the routes naturally adjust. However, if all visits can still be met, routes often require minimal or no adjustments to accommodate higher robustness or shorter time windows.

\begin{figure}[h] 
\centering \footnotesize
\label{tab:rev}
\begin{tabular}{cccccccccccccccccccccccc}
Wide && Narrow && Tight  \\ 
 \includegraphics[width=4.5cm]{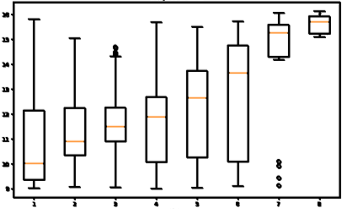} && \includegraphics[width=4.5cm]{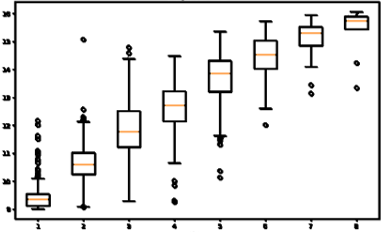} && \includegraphics[width=4.5cm]{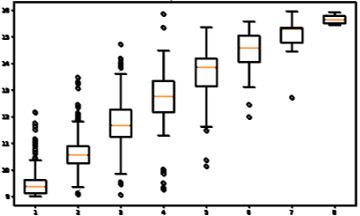}  \\
\end{tabular}
\caption{Daily schedules with medium robustness by time window size \textmd{(for Region T-20, Nursing)} \\}
\label{results_routes_schedules}
\end{figure}

We find that routes may remain similar across different robustness levels if the set of visited patients remains unchanged. The main adjustments typically occur in the schedules, with discrepancies becoming more noticeable as the day progresses. While dependent on the estimated deviations of service and travel times, the worst-case start times change on average as follows. Compared to no robustness, the worst-case start time for the 8th visit is 29.30 minutes later for low robustness and 33.74 minutes later for medium robustness. With high robustness, there was no 8th visit on any day. Furthermore, as depicted in Figure~\ref{results_routes_schedules}, decreasing time window lengths tend to increase the scheduled time between any two visits. 
In summary, scheduling buffers between visits enables shorter time windows and increased robustness, with limited profit loss as long as revenue remains unaffected.

\section{Summary and Conclusions}
\label{conclusion}

This work introduces a comprehensive decision tool for patient selection, caregiver assignment, and multi-day routing and scheduling, taking into account polyhedral uncertainty sets for travel and service times. It incorporates a variety of optimization criteria, including patient time windows, number of required visits, caregiver availability, compatibility, continuity of care, and profit. The tool is designed to ensure optimal treatment for all selected patients, considering the factors influenced by operational planning. In particular, all plans are robust, ensuring timely availability for patients and caregivers even in worst-case scenarios, which can be adjusted through uncertainty budgets. At the same time, profit is maximized, calculated as revenue minus wage cost and travel~cost.

We efficiently decompose the robust home healthcare routing and scheduling problem in several layers for a nested branch-and-price procedure. In the first layer, the problem is decomposed by caregiver, yielding weekly plans. Next, the problem of finding a weekly plan for each caregiver is decomposed by day, yielding daily routes for each. 
These  routes are computed using a branch-and-bound procedure, with traveling salesman subproblems (R-TSPTW), and a multi-knapsack problem guiding the acceptance decision.

Our case study with real data from a large US home healthcare agency demonstrates that, while runtimes increase with shorter time windows and the level of robustness, all instances can be optimally solved within reasonable computation times. Surprisingly, geographic factors seem secondary to revenue and feasibility in the acceptance and assignment of new patients. Compared to parameter settings absent of robustness and time windows, introducing a medium level of robustness and one-hour time windows, and thus smoother operations and a higher service level, leads to an average profit loss of 14.20\%. This is mainly due to decreased patient loads per caregiver per day. 
Nevertheless, some daily schedules can be robustified by scheduling more buffer times without affecting revenue.

Our work also points to several directions for future research. Firstly, the problem setting could be expanded. Alternative objectives such as preference satisfaction and workload balance, additional constraints for lunch breaks, caregiver dependencies, as well as uncertain caregiver absence are also relevant from a practical standpoint.
Secondly, runtimes could be further accelerated. While we solved the R-TSPTW subproblems using Gurobi, it seems promising, particularly for settings with high uncertainty, to employ a problem-specific solution method.
Thirdly, it would be worthwhile to investigate a variety of acceptance and assignment policies. This could include a comparison of our approach with advanced clustering approaches. More research is also needed on fair acceptance procedures for new patients.

\ACKNOWLEDGMENT{%
Funded by the Deutsche Forschungsgemeinschaft (DFG, German Research Foundation)
-- 277991500/GRK2201.
}

\addcontentsline{toc}{section}{References}
\bibliographystyle{template/informs2014} 
\bibliography{references}

\newpage

\begin{APPENDIX}{}

\section{Adversarial Problems: Dualization Scheme} \label{dual_scheme}
To apply the classical dualization scheme as proposed by \cite{bertsimas2004price}, we reformulate the service start times. In particular, we replace the binary variable $x$ with a new binary variable $z_{j'j}^{k, d, i} \in \lbrace0,1\rbrace$, which additionally indicates the sequence of patients treated by caregiver $k$ on day $d$. Specifically, $z_{j'j}^{kdi}$ takes a value of 1 if and only if the $i$th patient treated by caregiver $k$ on day $d$ is patient $j$, and the patient treated before $j$ was patient $j'$. Notably, we keep both patient indices due to the consideration of travel times. Now, we shift the focus from the start time of patient $i$ to the start time of the $i$th patient treated by caregiver $k$ on day $d$.

\begin{lemma}[Start Time Reformulation] \label{starttimes} Let $\bar{p}_0 := 0$ and $\hat{p}_0 := 0$. For a fixed caregiver and day, consider the binary variable $z_{j',j}^{i}$, where $z_{j',j}^{i}$ equals 1 if the $i$th visit is with patient $j \in \mathcal{I}$ and patient $j' \in \mathcal{I} \cup \lbrace 0 \rbrace$ was visited immediately before, and 0 otherwise. Then, the service start times can be computed as follows. 
\begin{enumerate}
    \item[a.] Compute the service start time $s_i$ for patient $i \in \mathcal{I}$, initializing $s_0$ as $\underaccent{\bar}{W}$, as follows: \begin{align} \label{start_formulation}     \stepcounter{alignnumber}     \setcounter{equation}{0}  s_i\geq \max \left \lbrace \underaccent{\bar}{T}_{{ i} }, \max_{j \in \mathcal{I} \cup \lbrace 0 \rbrace} \left(s_j + \bar{p}_i + \gamma_i^p \hat{p}_i + \bar{t}_{ij} + \gamma_{ij}^t \hat{t}_{ij} \right) x_{ji} \right\rbrace \end{align}
    \item[b.] 
   Compute the service start time $\tilde{s}_i$ for the $i$th visited patient, with $\underaccent{\bar}{T}_0 := \underaccent{\bar}{W}$, as follows: \begin{align}   \label{start_reformulation}       \setcounter{equation}{1}  \tilde{s}_i \geq \max \left \lbrace \sum_{j \in \mathcal{I}} \sum_{\substack{j' \in \mathcal{I} \\ \cup \lbrace 0 \rbrace}} z_{j'j}^{i}\underset{\bar{}}{T}_{j}, \max_{i' \leq i} \left \lbrace \sum_{j \in \mathcal{I}} \sum_{\substack{j' \in \mathcal{I} \\ \cup \lbrace 0 \rbrace}} z_{j'j}^{i'}\underset{\bar{}}{T}_{j'} + \sum_{l = i'}^i \sum_{j \in \mathcal{I}} \sum_{\substack{j' \in \mathcal{I} \\ \cup \lbrace 0 \rbrace}} (\bar{p}_j + \gamma_j^p\hat{p}_j + \bar{t}_{jj'} + \gamma_{jj'}^t \hat{t}_{jj'})z_{j'j}^{l} \right \rbrace \right \rbrace \end{align}  
\end{enumerate}
\end{lemma}
\vspace{0.5cm}

\begin{proof}{Proof of Lemma \ref{starttimes}.}
Let $\tilde{s}_0 := \underaccent{\bar}{W}$.
We prove by induction that (\ref{start_reformulation}) is equivalent to the recursive formulation of start times 
\begin{align}   \label{start_reformulation_rec}       \setcounter{equation}{2}  \tilde{s}_i \geq \max \left \lbrace \sum_{j \in \mathcal{I}} \sum_{j' \in \mathcal{I} \cup \lbrace 0 \rbrace} z_{j'j}^{i}\underset{\bar{}}{T}_{j}, \tilde{s}_{i-1} + \sum_{j \in \mathcal{I}} \sum_{j' \in \mathcal{I} \cup \lbrace 0 \rbrace} (\bar{p}_j + \gamma_j^p\hat{p}_j + \bar{t}_{jj'} + \gamma_{jj'}^t \hat{t}_{jj'})z_{j'j}^{i} \right \rbrace.  \end{align}

For $i = 1$, the terms are equivalent. Suppose the claim is true for all visits $1, \ldots, i-1$. Then replacing $\tilde{s}_{i-1}$ in (\ref{start_reformulation_rec}) with the formulation given by (\ref{start_reformulation}) shows the equivalence. It is evident that computing the service start times using (\ref{start_reformulation_rec}) is equivalent to using (\ref{start_formulation}).
\qed
\end{proof}
\vspace{0.5cm}

Let $n_{max}$ denote the maximum number of patients treated by any caregiver on any day. Using the binary variable $z$ and defining $p_0 := 0$, $\hat{p}_0 := 0$, $\bar{p}_0 := 0$, we can reformulate constraints~{\ref{model1con1}} as:

\begin{align}     \stepcounter{alignnumber}     \setcounter{equation}{0} 
    s_{i^*}^{kd} & \leq \sum_{j=1}^n \sum_{j'=1}^n  \bar{T}_{jd} z_{j'j}^{k,d,{i^*}}  & \forall k \in \mathcal{K}\, \forall d \in \mathcal{D}\,\forall i^* \in [n_{max}]  \\
    s_{i^*}^{kd} & \geq \underaccent{\bar}{W}_{kd}  + \max_{\substack{p \in \mathcal{U}^p \\t \in \mathcal{U}^t}} \left\lbrace \sum_{i = 1}^{i^*}\sum_{j=1}^n \sum_{j'=0}^n (p_j + t_{j'j})z_{j'j}^{k,d,i} \right \rbrace & \forall k \in \mathcal{K}\, \forall d \in \mathcal{D}\, \forall i^* \in [n_{max}] \label{max1} \\
    s_{i^*}^{kd} & \geq \sum_{j=1}^n \sum_{j'=0}^n  \underaccent{\bar}{T}_{jd} z_{j'j}^{k,d,i'}   + \max_{\substack{p \in \mathcal{U}^p \\t \in \mathcal{U}^t}} \left\lbrace \sum_{i = i'}^{i^*}\sum_{j=1}^n \sum_{j'=0}^n (p_j + t_{j'j})z_{j'j}^{k,d,i} \right\rbrace & \forall k \in \mathcal{K}\, \forall d \in \mathcal{D}\,\forall i^* \in [n_{max}] \forall i' \in [i^*]  \label{max2}
\end{align}

Now, replacing the maximization problems in \ref{max1} and \ref{max2} by their LP duals, we can replace constraints~\hyperref[R-HHCRSP]{{\ref{model1con1}}} by the following linear formulation:

\begin{align*}
    s_{i^*}^{kd} & \leq \sum_{j=1}^n \sum_{j'=1}^n  \bar{T}_{jd} z_{j'j}^{k,d,{i^*}}    & \forall k \in \mathcal{K}\, \forall d \in \mathcal{D}\,\forall i^* \in [n_{max}]  \\
    s_{i^*}^{kd} & \geq \underaccent{\bar}{W}_{kd}  + \sum_{i = 1}^{i^*}\sum_{j=1}^n \sum_{j'=0}^n (\bar{p}_j + \bar{t}_{j'j})z_{j'j}^{k,d,i}  \nonumber \\& \multicolumn{2}{l}{$\quad \,\,\,\,\, \quad + \Gamma^p y_p^{(k,d,i^*,0)} + \Gamma^t y_t^{(k,d,i^*,0)} + \sum_{j=1}^n (y^p)_{j}^{^{(k,d,i^*,0)}} + \sum_{j=1}^n \sum_{j'=0}^n  (y^t)_{j'j}^{^{(k,d,i^*,0)}} $} \nonumber \\&& \forall k \in \mathcal{K}\, \forall d \in \mathcal{D}\,\forall i^* \in [n_{max}] \\
    y_p^{(k,d,i^*,0)} + (y^p)_{j}^{^{(k,d,i^*,0)}} & \geq \sum_{i = 1}^{i^*} \sum_{j'=1}^n \hat{p}_j z_{j'j}^{k,d,i} & \forall k \in \mathcal{K}\, \forall d \in \mathcal{D}\,\forall i^* \in [n_{max}]\forall j \in N \\
    y_t^{(k,d,i^*,0)} + (y^p)_{j'j}^{^{(k,d,i^*,0)}} & \geq \sum_{i = 1}^{i^*} \hat{t}_{j'j} z_{j'j}^{k,d,i}   & \forall k \in \mathcal{K}\, \forall d \in \mathcal{D}\,\forall i^* \in [n_{max}] \forall j' \in N_0 \forall j \in N \\
    s_{i^*}^{kd} & \multicolumn{2}{l}{$\geq \sum_{j=1}^n \sum_{j'=0}^n  \underaccent{\bar}{T}_{jd} z_{j'j}^{k,d,i'}  + \sum_{i = i'}^{i^*}\sum_{j=1}^n \sum_{j'=0}^n (\bar{p}_j + \bar{t}_{j'j})z_{j'j}^{k,d,i}$} \nonumber \\& \multicolumn{2}{l}{$\qquad + \Gamma^p y_p^{(k,d,i^*,i')} + \Gamma^t y_t^{(k,d,i^*,i')}  + \sum_{j=1}^n (y^p)_{j}^{^{(k,d,i^*,i')}} + \sum_{j=1}^n \sum_{j'=0}^n  (y^t)_{j'j}^{^{(k,d,i^*,i')}}$}   \nonumber \\&&  \forall k \in \mathcal{K}\, \forall d \in \mathcal{D}\,\forall i^* \in [n_{max}] \forall i' \in [i^*] \\
    y_p^{(k,d,i^*,i')} + (y^p)_{j}^{^{(k,d,i^*,i')}} & \geq \sum_{i = i'}^{i^*} \sum_{j'=1}^n \hat{p}_j z_{j'j}^{k,d,i}  & \forall k \in \mathcal{K}\, \forall d \in \mathcal{D}\,\forall i^* \in [n_{max}] \forall i' \in [i^*] \forall j \in N \\
    y_t^{(k,d,i^*,i')} + (y^p)_{j'j}^{^{(k,d,i^*,i')}} & \geq \sum_{i = i'}^{i^*} \hat{t}_{j'j} z_{j'j}^{k,d,i}  & \forall k \in \mathcal{K}\, \forall d \in \mathcal{D}\,\forall i^* \in [n_{max}] \forall i' \in [i^*] \forall j' \in N_0 \forall j \in N \\
    y &\geq 0 
\end{align*}

This formulation requires $|\mathcal{I}|^2\cdot|\mathcal{K}|\cdot|\mathcal{D}|\cdot n_{max}$ binary variables, $\mathcal{O}\left(|\mathcal{I}|^2\cdot|\mathcal{K}|\cdot|\mathcal{D}|\cdot \left(n_{max}\right)^2 \right)$ new continuous variables and replaces constraints~{\ref{model1con1}} with $\mathcal{O}\left(|\mathcal{I}|^2\cdot|\mathcal{K}|\cdot|\mathcal{D}|\cdot \left(n_{max}\right)^2 \right)$ new constraints.

\section{Nested Branch-and-Price Procedure}
In the following, we present the LP duals for the linear relaxation of both the first- and second-level master problems, along with the corresponding pricing problems.

\subsection{LP Duals} \label{appendix:duals}

The LP dual of the linear relaxation of the first-level master problem (\ref{MP-1}) -- (\ref{MP1-last}) can be stated as follows: 

\begin{align}     \stepcounter{alignnumber}     \setcounter{equation}{0} 
    \min\, & \sum_{k\in\mathcal{K}} v_k +\sum_{i\in\mathcal{I}_R} w_i \label{dualobj}  \\
    \text{s.t.}\, & v_k +  \sum_{i\in\mathcal{I}_R} a_{k,P}^i \, w_i \geq Rev_{k,P} - C(T)_{k,P} - C(W)_{k,P}    & \forall k \in \mathcal{K}\forall P \in \mathcal{P}_k \\
    & w \geq 0 \label{duallast}
\end{align}

The dual perspective provides valuable insights into the problem. It aims to find an optimal value for each caregiver and each new patient. For each weekly plan, the total value generated by the caregiver and the newly accepted patients he or she visits must be at least as high as the profit generated by the weekly plan. The LP dual aims to minimize the aggregate value associated with all caregivers and new patients.

The LP dual of the linear relaxation of the second-level master problem (\ref{MP2-first}) -- (\ref{MP2-last})  can be stated as follows: 

\begin{align}     \stepcounter{alignnumber}     \setcounter{equation}{0} 
    \min\, & \sum_{d\in\mathcal{D}} u_d +  \sum_{i\in\mathcal{I}^{\mathcal{R}}} v_i z_i + \sum_{i\in \mathcal{I}^{\mathcal{R}}}\sum_{\substack{\tilde{d} \in [|\mathcal{D}| - d_i]_0 \cap [d_i -d]_0}}  q_{i\tilde{d}} \label{dual22}  \\
    \text{s.t.}\, &  u_d +    \sum_{i\in\mathcal{I}^{\mathcal{R}}} a_{d, R}^iz_i + \sum_{i\in\mathcal{I}^{\mathcal{R}} }\sum_{\substack{\tilde{d} \in [|\mathcal{D}| - d_i]_0 \cap [d_i -d]_0}}  a_{d, R}^i q_{i\tilde{d}} \label{formula_days}\\
    & \qquad - \sum_{i \in \mathcal{I}^{\mathcal{R}}} \sum_{\tilde{d} \in \mathcal{D}} a_{d, R}^i y_{i\tilde{d}} + \sum_{i \in \mathcal{I}^{\mathcal{R}}} v_i a_{d, R}^i y_{id}  \\
     & \qquad \geq Rev_{d,R} - C(T)_{d,R} - C(W)_{d,R} - \sum_{i \in \mathcal{I}^{\mathcal{R}}} \frac{w_{i}^*}{v_i} a_{d,R}^i & \forall d \in \mathcal{D}\, \forall R \in \mathcal{R}_d \\
    &  z, y, q \geq 0 & \label{dual22end}
\end{align}

The second-level dual LP evaluates the constraints of a single route per day ($u$), the fulfillment of the exact number of required visits for each new patient ($z, y$), and the adherence to the minimum required time period between two consecutive visits ($q$).

\subsection{Pricing problem} 
\label{1stlevel_pricing_MIP} 

\label{2ndlevel_pricing_MIP}
The second-level pricing problem is solved independently for each fixed caregiver $k \in \mathcal{K}$ and day $d \in \mathcal{D}$. It generates a feasible daily route that maximizes the profit, minus the value of the newly accepted patients (determined by the first- and second-level dual values).

\begin{align}     \stepcounter{alignnumber}     \setcounter{equation}{0}  \label{pricing_pricing}
    \max\, &   \sum\limits_{\substack{i \in \mathcal{I}_{k}}} R_{ki}\delta_{i} -     \sum\limits_{\substack{i, j \in \mathcal{I}_{k}}} C^T_{dij} x_{ij} -   C^W_k    \left(s_{n+1} -\underaccent{\bar}{W}_{{ kd} } \right) 
          \\ \nonumber 
       & \multicolumn{2}{l}{ $ - \left( \sum_{i \in \mathcal{I}^{\mathcal{R}}} \frac{w_{i}^*}{v_i} + \sum_{i\in\mathcal{I}^{\mathcal{R}}} z_i^* + \sum_{i \in \mathcal{I}^{\mathcal{R}}} \left( \sum_{\tilde{d} \in \mathcal{D}} y_{i\tilde{d}}^* -  v_i y_{id}^*  \right)
 + \sum_{i\in\mathcal{I}_{k}}\sum_{\substack{\mu\in [|\mathcal{D}| - d_i]: \\ \mu \in [d_i -d]_0}} q_{i,\mu}^*\, \right) \delta_i$} \\ \nonumber \\
     \text{s.t.} \, 
    & 
    s_{i,\gamma^p +   \zeta^p, \gamma^t +   \zeta^t}    \geq  s_{j,\gamma^p , \gamma^t}  + (\bar{p}_j +   \zeta^p\hat{p}_j + \bar{t}_{ji} +   \zeta^t\hat{t}_{ji})  x_{ji}   \\ \nonumber 
    &\multicolumn{2}{l}{ $\qquad\qquad \qquad + ( \underaccent{\bar}{T}_{{ id} } - \bar{T}_{{ jd} } )(1-  x_{ji} ) \qquad\qquad\quad\quad\quad\quad\quad \,\,\, \forall i,j\in\mathcal{I}_{k}  \, \forall \gamma^p \in [\Gamma^p]_0 \, \forall \gamma^t \in [\Gamma^t]_0$ } \nonumber \\ \nonumber 
    & \multicolumn{2}{l}{$\qquad\qquad\qquad\quad\qquad\qquad\,\,\qquad\qquad\qquad\qquad \,\,\,\,\,\, \forall   \zeta^p,   \zeta^t \in \lbrace 0, 1 \rbrace: \gamma^p +   \zeta^p \leq \Gamma^p, \gamma^t +   \zeta^t \leq \Gamma^t$}  \\
     &  s_{i,0,   \zeta^t}   \geq  \left(\underaccent{\bar}{W}_{{ kd} } + \bar{t}_{0i} +   \zeta^t\hat{t}_{0i}\right)  x_{0i}  + (\underaccent{\bar}{W}_{{ kd} } - \bar{T}_{{ id} })(1- x_{i,n+1} )   & \forall i\in\mathcal{I}_{k}\,  \forall  \zeta^p \zeta^t \in \lbrace 0, 1 \rbrace  \\
            &  s_{i,\gamma^p, \gamma^t}  \in [ \underaccent{\bar}{T}_{{ id} }, \bar{T}_{{ id} }] & \forall i \in \mathcal{I}_{k} \, \forall \gamma^p \in [\Gamma^p]_0 \, \forall \gamma^t \in [\Gamma^t]_0 \\
              &  s_{n+1}   \geq  s_{i,\gamma^p - \zeta^p, \Gamma^t - \zeta^t}  + (\bar{p}_i + \zeta^p\hat{p}_i + \bar{t}_{i,n+1} + \zeta^t \hat{t}_{i,n+1})  x_{i,n+1}  &  \forall i\in\mathcal{I}_{k}\,  \forall   \zeta^t \in \lbrace 0, 1 \rbrace  \\ 
     & s_{n+1}  \in [\underaccent{\bar}{W}_{{ kd} } , \bar{W}_{{ kd} }] &  \\
    &  \sum_{j\in\mathcal{I}_{k}} x_{ij}  = \delta_i &   \forall i \in \mathcal{I}_{k} \\[3pt]
                                  &  \sum_{j\in \mathcal{I}_{k}} x_{ij}  \leq C_i^{kd}  &    \, \forall i\in \mathcal{I}_{k}\\[3pt]
    & \sum_{j\in\mathcal{I}_{k}} x_{ji}   -\sum_{j\in\mathcal{I}_{k}} x_{ij}  = 0, \sum_{i\in\mathcal{I}_{k}} x_{ii}  = 0 &    \, \forall i\in \mathcal{I}_{k} \\[3pt]
    & \sum_{\substack{i \in \mathcal{I}_{k} \\ \cup \lbrace n+1 \rbrace}} x_{0,i}  = 1, \sum_{\substack{i\in \mathcal{I}_{k} \\ \cup \lbrace0 \rbrace}} x_{i,n+1}  = 1 &     \\[3pt]
   & x_{ij} , \delta_i \in \lbrace 0,1 \rbrace &    \, \forall i, j\in \mathcal{I}_{k} \cup\lbrace 0 \rbrace\cup \lbrace n+1 \rbrace \\
   & s_{i,\gamma^p, \gamma^t} , s_{n+1} \geq 0 &  \, \forall i \in \mathcal{I}_{k} \, \forall \gamma^p \in [\Gamma^p]_0 \, \forall \gamma^t \in [\Gamma^t]_0
\end{align}

The second-level pricing problem is a MILP with $\mathcal{O}\left(n^2  \right) $ binary variables, $\mathcal{O} \left( n  \cdot|\Gamma^p|\cdot|\Gamma^t| \right) $ continuous variables and $\mathcal{O}(n^2  \cdot |\Gamma^p|\cdot|\Gamma^t|)$ constraints.

\section{Case Study} \label{compexp}
Appendix~\ref{data_description} provides detail on the data and Appendix~\ref{results_description} provides additional details on the computational results.

\subsection{Data.} \label{data_description}
The original data set encompasses detailed information on 166,641 visits made by 807 caregivers to 7,405 patients across eight regional branches and five different disciplines. Data was collected by a US national healthcare agency headquartered in Dallas, Texas over a 36-week period from July 21st, 2018 to March 18th, 2019. Figure~\ref{fig:timelycourse} illustrates the weekly visit counts for the entire data set. To ensure data completeness, we removed weeks with less than 4,000 visits, resulting in 23 considered weeks. The data comprises information on the patient locations, required services, and timely availability for each patient, as well as the qualification, employment status, and compensation type for each caregiver. Table~\ref{tab:caregivers} summarizes the type of employment and compensation for all caregivers in the data set. 

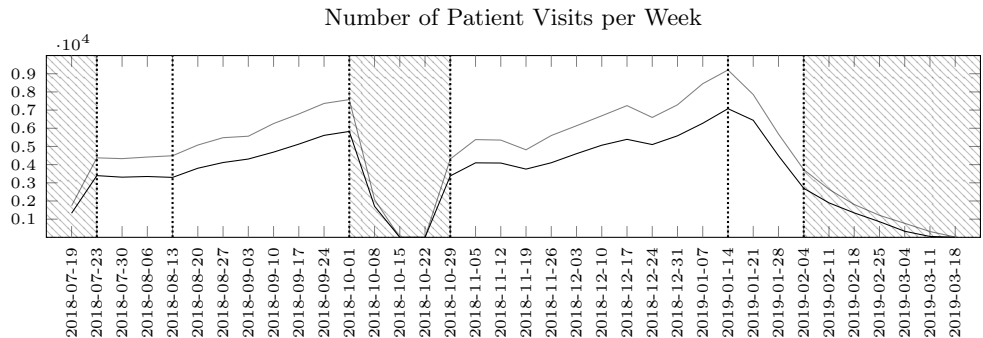
\begin{figure}[b]
    \centering
    \begin{tikzpicture}
\begin{axis}[width=14cm,height=4cm,
    title={\footnotesize Number of Patient Visits per Week},
    ticklabel style = {font=\tiny}, 
    xmin=0, xmax=37,
    ymin=0, ymax=10000,
     xtick=\empty,
    extra x ticks={1,	2,	3,	4,	5,	6,	7,	8,	9,	10,	11,	12,	13,	14,	15,	16,	17,	18,	19,	20,	21,	22,	23,	24,	25,	26,	27,	28,	29,	30,	31,	32,	33,	34,	35,	36 },
    extra x tick labels={2018-07-19,	2018-07-23,	2018-07-30,	2018-08-06,	2018-08-13,	2018-08-20,	2018-08-27,	2018-09-03,	2018-09-10,	2018-09-17,	2018-09-24,	2018-10-01,	2018-10-08,	2018-10-15,	2018-10-22,	2018-10-29,	2018-11-05,	2018-11-12,	2018-11-19,	2018-11-26,	2018-12-03,	2018-12-10,	2018-12-17,	2018-12-24,	2018-12-31,	2019-01-07,	2019-01-14,	2019-01-21,	2019-01-28,	2019-02-04,	2019-02-11,	2019-02-18,	2019-02-25,	2019-03-04,	2019-03-11,	2019-03-18},
    extra x tick style={
           tick label style={rotate=90} 
           },
    ytick={1000, 2000, 3000, 4000, 5000, 6000, 7000, 8000, 9000},
    legend pos=north west,
    grid style=dashed,
]

\addplot[color = gray,
    ]
    coordinates {
    (1,1736)	(2,4373)	(3,	4330)	(4,4417)	(5,4489)	(6,5077)	(7,5482)	(8,5563)	(9,6259)	(10,6786)	(11,7365)	(12,7578)	(13,	2080)	(14,0)	(15,0)	(16,4309)	(17,5377)	(18,5350)	(19,4816)	(20,5601)	(21,6139)	(22,6686)	(23,7246)	(24,6595)	(25,7291)	(26,8453)	(27,9206)	(28,7856)	(29,5686)	(30,3727)	(31,2643)	(32,1807)	(33,1210)	(34,771)	(35,319)	(36,18)
    };
    
\addplot[
]
coordinates {
(1,1327)	(2,3395)	(3,	3310)	(4,3347)	(5,3297)	(6,3801)	(7,4115)	(8,4308)	(9,4688)	(10,5127)	(11,5609)	(12,5826)	(13,	1709)	(14,0)	(15,0)	(16,3380)	(17,4098)	(18,4087)	(19,3754)	(20,4103)	(21,4602)	(22,5067)	(23,5391)	(24,5104)	(25,5585)	(26,6275)	(27,7078)	(28,6439)	(29,4487)	(30,2697)	(31,1897)	(32,1345)	(33,861)	(34,334)	(35,48)	(36,3)
};
\addplot[thick, samples=50, densely dotted,domain=0:6,black, name path=three] coordinates {(2,0)(2,10000)};
\addplot[thick, samples=50, densely dotted,domain=0:6,black, name path=three] coordinates {(5,0)(5,10000)};
\addplot[thick, samples=50, densely dotted,domain=0:6,black, name path=three] coordinates {(12,0)(12,10000)};
\addplot[thick, samples=50, densely dotted,domain=0:6,black, name path=three] coordinates {(16,0)(16,10000)};
\addplot[thick, samples=50, densely dotted,domain=0:6,black, name path=three] coordinates {(27,0)(27,10000)};
\addplot[thick, samples=50, densely dotted,domain=0:6,black, name path=three] coordinates {(30,0)(30,10000)};
\addplot [pattern=north west lines,densely dotted,opacity=0.5, draw = black] coordinates {
         (0,\pgfkeysvalueof{/pgfplots/ymax})
         (2,\pgfkeysvalueof{/pgfplots/ymax})
        }
        \closedcycle;
\addplot [pattern=north west lines,densely dotted,opacity=0.5, draw = black] coordinates {
         (30,\pgfkeysvalueof{/pgfplots/ymax})
         (37,\pgfkeysvalueof{/pgfplots/ymax})
        }
        \closedcycle;     
\addplot [pattern=north west lines,densely dotted,opacity=0.5, draw = black] coordinates {
         (12,\pgfkeysvalueof{/pgfplots/ymax})
         (16,\pgfkeysvalueof{/pgfplots/ymax})
        }
        \closedcycle;     
\end{axis}

\end{tikzpicture}
    \caption{Temporal trends in weekly patient visits: Overall visits depicted in gray, regular visits in black; shaded regions indicate periods with fewer than 4,000 visits.}
    \label{fig:timelycourse}
\end{figure}

 \begin{table}[h]
  \caption{Caregivers}
    \centering \small
    \begin{tabular}{rrrrrrrrrrrrrrrrrrrrrrrr} \toprule
         \multicolumn{7}{l}{\bfseries Employment} && \multicolumn{7}{l}{\bfseries Compensation} \\ \cmidrule{1-7} \cmidrule{9-15}
        full time && PRN && part time && other contract && per visit && hourly && salaried && other contract \\\cmidrule{1 - 1}  \cmidrule{3-3} \cmidrule{5-5} \cmidrule{7-7} \cmidrule{9-9} \cmidrule{11-11}  \cmidrule{13-13}  \cmidrule{15-15}
          51.93\% && 18.98\% && 2.49\% && 26.60\% && 37.41\% && 23.57\% && 12.68\% && 26.34\%  \\ \bottomrule
    \end{tabular}
    \label{tab:caregivers}
\end{table}

\begin{figure}[b]
    \centering
    \includegraphics[scale = 0.45]{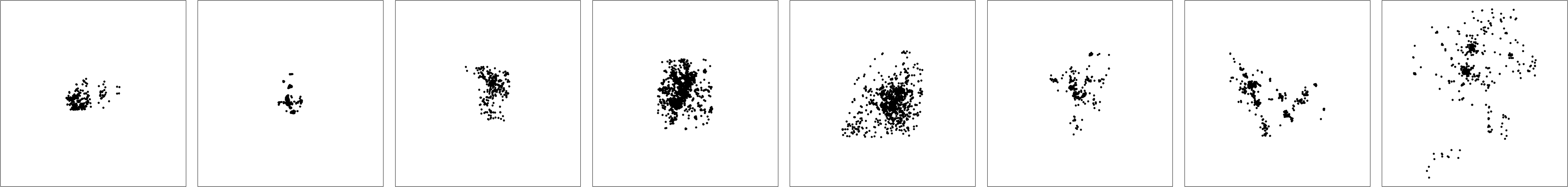}
    \caption{Scatter plots of patient locations by regional branch \qquad \qquad \qquad \qquad \qquad \qquad \qquad \qquad \qquad \qquad \qquad \qquad \qquad \qquad \qquad  (arranged \textmd{(arranged in increasing order with respect to average travel time)}}
    \label{fig:regions}
\end{figure}

\paragraph{Regional branches.} The original data set comprises visits from eight regional branches located in four US states. With the exception of seven caregivers, all others exclusively visited patients in a single region. A summary of the eight regions is presented in Table~\ref{tab:regions}, and a scatter plot of patient locations within each region is shown in Figure~\ref{fig:regions}. Variability in travel distances and times within each region influences the home healthcare routing and scheduling problem. Accordingly, we categorized the regions into three groups: town-shorttravel, city-midtravel, and town-longtravel. Out of each group, we chose up to one candidate for our instance generation.

\begin{table}[t]
    \centering \small
    \caption{Regional branches}
    \begin{tabular}{crrrrrrrrrrrr} \toprule
            &&  Patients && { Travel time} \newline {(avg.)} && { Distance} \newline {(avg.)} &&  Population && Group \\  \cmidrule{3-3} \cmidrule{5-5} \cmidrule{7-7} \cmidrule{9-9} \cmidrule{11-11} 
        \bf Region-T19: && 729&& {19.55 min} && {21.91 km} && $<500,000$ && town-shorttravel   \\ 
        \bf Region-T20: && 1,138&& {20.22 min} && {22.60 km} && $< 100,000$ && town-shorttravel  \\ 
        \bf Region-T30: && 525&& {30.82 min} && {28.98 km} && $< 500,000$ && town-shorttravel  \\ 
        \bf Region-T33: && 1,856&& {33.94 min} && {38.03 km} && $< 1,000,000$ && city-midtravel  \\ 
        \bf Region-T36: && 1,677&& {36.06 min} && {41.87 km} && $< 1,500,000$ && city-midtravel \\ 
        \bf Region-T40: && 398&& {40.64 min} && {42.63 km} && $< 10,000$ && town-longtravel \\ 
        \bf Region-T48: && 559&& {48.06 min} && {54.18 km} && $< 100,000$ && town-longtravel \\ 
        \bf Region-T60: && 523&& {60.98 min} && {74.69 km} && $< 100,000$ && town-longtravel \\ \bottomrule
    \end{tabular}
    \label{tab:regions}
\end{table}

\begin{table}[b]
    \centering \small
    \caption{Services \textmd{(excluding 1.78\% categorized as 'Other')}}
    \begin{tabular}{crrrrrrrrrrrrrrrrrrrrrrrrr} \toprule
          && {Total visits}   && \multicolumn{3}{c}{Frequency} (per patient) &&  \multicolumn{3}{c}{Visit types} && {Change rate} && {Wage}     \\    \cmidrule{5-7}  \cmidrule{9-11}   
          && &&  Weekly &&   LoS && Routine && Other  &&   &&    \\  \cmidrule{3-3} \cmidrule{5-5} \cmidrule{7-7} \cmidrule{9-9}  \cmidrule{11-11}  \cmidrule{13-13}  \cmidrule{15-15} 
             \bf SN: && 41,889 &&  1 visit(s)  && 5 weeks && 100\,\%  && 0\,\% &&   0.65 && \$26.86    \\
            \bf RN: && 18,361 && 1 visit(s) &&  2 weeks && 18\,\%  && 82\% && 4.84   && \$42.80    \\ 
            \bf LPN: && 6,959 &&  1 visit(s) && 4 weeks && 98\,\%  &&  2\,\%  &&  0.75  &&  \$26.86  \\ 
         \bf PSY: && 6,446  && 1 visit(s) &&   9 weeks&& 80\,\%   && 20\,\%  && 0.18 && \$59.94  \\  \cmidrule{3-3} \cmidrule{5-5} \cmidrule{7-7} \cmidrule{9-9}  \cmidrule{11-11}  \cmidrule{13-13}  \cmidrule{15-15} 
        \bf PTA: && 34,639 && 2 visit(s) &&  5 weeks && 100\,\%  && 0\,\%&& 0.31   && \$31.01      \\ 
           \bf PT: && 25,780  &&  1 visit(s) &&  3 weeks&& 46\,\%  && 54\,\%  && 1.36  && \$47.10  \\ \cmidrule{3-3} \cmidrule{5-5} \cmidrule{7-7} \cmidrule{9-9}  \cmidrule{11-11}  \cmidrule{13-13}  \cmidrule{15-15} 
          \bf OT: && 12,037  && 1 visit(s)  && 2 weeks&& 49\,\% && 51\,\% &&   1.00 &&  \$44.61  \\
           \bf COTA: && 8,869   &&   1.5 visit(s) && 3 weeks && 100\,\%  && 0\,\%  &&  0.43 && \$20.71   \\ \cmidrule{3-3} \cmidrule{5-5} \cmidrule{7-7} \cmidrule{9-9}  \cmidrule{11-11}  \cmidrule{13-13}  \cmidrule{15-15} 
     \bf ST:  && 6,195 &&   1 visit(s)   &&4 weeks&& 71\,\%  && 29\,\% &&  0.28  && \$38.01     \\ 
        \bf HHA: && 2,502 && 1 visit(s) &&  3 weeks&& 100\,\% &&  0\,\% &&  0.36  && \$14.15   \\ 
     \bf MSW:    && 1,074 && 1 visit(s) && 1 weeks && 15\,\%  &&  85\,\% && 9.07  && \$30.17 \\ 
         \bottomrule
    \end{tabular}
    \label{tab:services}
\end{table}

\paragraph{Services.} Patient visits fall into the following disciplines: nursing (45.33\%  - comprised of skilled nursing [SN], registered nursing [RN], licensed practical nursing [LPN], and psychiatric nursing [PSY]), physical therapy (36.26\% - comprised of physical therapy [PT] and physical therapy assistance [PTA]), occupational therapy (12.55\% - comprised of occupational therapy  [OT] and occupational therapy aid [COTA]), and other categories (5.87\% - including speech therapy [ST], home health aid [HHA], and medical social work [MSW]). Qualifications required for the different services vary. For example, PSY services necessitate PSY-qualified nurses, RN services require RN-qualified nurses, while LPN and SN services can be delivered by both LPN- and RN-qualified nurses. The qualification hierarchy is denoted as $RN \succcurlyeq LPN$, $RN \succcurlyeq SN$, $PT\succcurlyeq PTA $, and $OT\succcurlyeq COTA$. The visits are further classified into various types, including routine visits (75.91\%), reassessment visits (9.88\%), initial visits (7.98\%), discharge visits (2.62\%), and other types of visits (3.62\%). In our computational study, we focus exclusively on scheduling routine visits. 
Tables~\ref{tab:services} summarizes key characteristics of each service that impact the home healthcare routing and scheduling problem, such as mean weekly visits and length of stay (LoS) per patient, visit type distribution, change rate (average ratio of new patients to existing patients), wage (according to \cite{us20220occupational}), and expected service time.

\subsection{Results} \label{results_description}
In addition to the results presented in the main part, Table~\ref{tab:results_runtime} summarizes the runtimes for all instances, averaged over three weeks. Similarly, Table~\ref{tab:results_summary} summarizes the profit across all instances.
\begin{table}[h]
\centering \footnotesize
\caption{Average runtimes by region, discipline, level of robustness, and time window size \qquad \qquad \qquad \qquad \qquad \qquad \qquad \qquad \qquad \qquad \qquad \qquad \qquad \qquad \qquad 
 \textmd{(in seconds)}}
\label{tab:results_runtime}
\begin{tabular}{rrlrrrrrrrrrrrrrrrrrrrrrrrrrrrrrrrrrrr}
\toprule
     &&  && \multicolumn{5}{c}{  Physical Therapy}  && \multicolumn{5}{c}{ Nursing}           \\   \cmidrule{5 - 9} \cmidrule{11 - 15} 
      &&  && Wide  && Narrow && Tight && Wide && Narrow && Tight        \\
\cmidrule{1 - 1} \cmidrule{3 - 3}\cmidrule{5 - 5}\cmidrule{7 - 7}\cmidrule{9 - 9}\cmidrule{11 - 11} \cmidrule{13 - 13} \cmidrule{15-15} 
  Region T-20   && None && 10.50 && 3.83 && 3.99 && 39.51 && 80.57 && 36.69 \\  
   && Low   && 177.67 && 12.67 && 15.20 && 411.04 && 13.86 && 42.43 \\
   && Medium   && 1,116.91 && 205.52 && 37.06 && 6,758.96 && 11.60 && 434.74 \\ 
   &&  High  && 25.67 && 3.45 && 3.43 && 51.86 && 187.55 && 11.61 \\
 \cmidrule{1-1} \cmidrule{3-3} \cmidrule{5 - 9} \cmidrule{11 - 15} 
 Region T-60   && None  && 9.92 && 1.33 && 2.91 && 21.00 && 1.39 && 1.52 \\
  && Low   && 81.65 && 5.45 && 9.54 && 159.62 && 10.16 && 11.41 \\
   && Medium   && 260.32 && 14.03 && 17.59 && 373.75 && 30.28 && 30.34 \\
 && High   && 4.35 && 1.27 && 2.25 && 4.67 && 1.24 && 0.86 \\ 
\bottomrule
\end{tabular}
\end{table}
\begin{table}[h]
\centering \footnotesize
\caption{Average profit by region, discipline, level of robustness, and time window size \qquad \qquad \qquad \qquad \qquad \qquad \qquad \qquad \qquad \qquad \qquad \qquad \qquad \qquad \qquad \textmd{(in USD)} }
\begin{tabular}{rrlrrrrrrrrrrrrrrrrrrrrrrrrrrrrrrrrr}
\toprule
     &&  && \multicolumn{5}{c}{ Nursing} && \multicolumn{5}{c}{ Physical Therapy}           \\   \cmidrule{5 - 9} \cmidrule{11 - 15} 
      &&  && Wide && Narrow && Tight && Wide && Narrow && Tight        \\
\cmidrule{1 - 1} \cmidrule{3 - 3}\cmidrule{5 - 9}\cmidrule{11 -15}
Region T-20 && None && 17,634 && 17,872 && 17,625 && 6,939 && 6,841 && 6,841 \\
&& Low && 17,634 && 17,625 && 17,620 && 6,914 && 6,817 && 6,641 \\
&& Medium && 17,634 && 17,625 && 17,620 && 6,819 && 6,731 && 6,624 \\
&& High && 17,634 && 17,624 && 17,620 && 5,448 && 5,394 && 5,391 \\ 
\cmidrule{1-1} \cmidrule{3-3} \cmidrule{5-9} \cmidrule{11-15}
Region T-60 && None && 4,525 && 4,200 && 4,290 && 6,202 && 5,458 && 5,447 \\
&& Low && 4,015 && 3,803 && 3,325 && 5,859 && 5,054 && 4,944 \\
&& Medium && 3,723 && 3,329 && 3,138 && 5,361 && 4,946 && 4,944 \\
&& High && 3,458 && 3,231 && 3,138 && 5,041 && 4,946 && 4,944 \\
\bottomrule
\end{tabular}
\label{tab:results_summary}
\end{table}

\end{APPENDIX}

\end{document}